\documentclass[12pt]{article}%
\usepackage{graphicx}
\usepackage{amsmath}
\usepackage{amsthm}
\usepackage{latexsym}
\usepackage{amsfonts}
\usepackage{amssymb}%
\newtheorem{theorem}{Theorem}[section]
\newtheorem{proposition}[theorem]{Proposition}

\newtheorem{corollary}[theorem]{Corollary}
\newtheorem{remark}[theorem]{Remark}

\begin{document}

\title{Sharp asymptotics for the partition function\\of some continuous-time directed polymers}
\author{
Agnese Cadel 
\hspace{0.3cm} Samy Tindel\\
{\small \textit{Institut Elie Cartan, Universit\'e de Nancy 1}}\\
{\small \textit{BP 239, 54506-Vandoeuvre-l\`es-Nancy, France}}\\
{\small \texttt{[cadel,tindel]@iecn.u-nancy.fr}}\vspace*{0.05in}
\and Frederi Viens 
\thanks{This author's research partially supported by NSF grant no.: 0204999.}\\
{\small \textit{Dept. Statistics \& Dept. Mathematics, Purdue University}}\\
{\small \textit{150 N. University St., West Lafayette, IN 47907-2067, USA}}\\
{\small \texttt{viens@purdue.edu}}\vspace*{0.05in}
	   }

\maketitle

\begin{abstract}
This paper is concerned with two related types of directed polymers in a
random medium. The first one is a $d$-dimensional Brownian motion living in a
random environment which is Brownian in time and homogeneous in space. The
second is a continuous-time random walk on ${\mathbb{Z}}^{d}$, in a random
environment with similar properties as in continuous space, albeit defined
only on ${\mathbb{R}}_{+}\times{\mathbb{Z}}^{d}$. The case of a space-time
white noise environment can be acheived in this second setting. By means of
some Gaussian tools, we estimate the free energy of these models at low
temperature, and give some further information on the strong disorder regime
of the objects under consideration.
\end{abstract}

\vspace{1cm}



\vspace{0.5cm}


\section{Introduction}

\subsection{Background, models, and motivation}

Models for directed polymers in a random environment have been introduced in
the physical literature \cite{DSp,HH,IS,Pi} for two main reasons. First, they
provide a reasonably realistic model of a particle under the influence of a
random medium, for which a number of natural questions can be posed, in terms
of the asymptotic behavior for the path of the particle. The second point is
that, in spite of the fact that polymers seem to be some more complicated
objects than other disordered systems such as spin glasses, a lot more can be
said about their behavior in the low temperature regime, as pointed out in
\cite{FH,HH}. At a mathematical level, after two decades of efforts, a
substantial amount of information about different models of polymer is now
available, either in discrete or continuous space settings (see
\cite{CY2,Peter,RT} and \cite{CH,Me} respectively).

\vspace{0.3cm}

The current article can be seen as a part of this global project consisting in
describing precisely the polymer's asymptotic behavior, beyond the spin glass
case. Except for some toy models such as the {\small REM} or {\small GREM}
\cite{Bo,Tbk}, little is known about the low temperature behavior of the free
energy for spin glasses systems, at least at a completely rigorous level. We
shall see in this paper that polymer models are amenable to computations in
this direction: we work to obtain some sharp estimates on the free energy of
two different kind of polymers in continuous time, for which some scaling
arguments seem to bring more information than in the discrete time setting.
Here, in a strict polymer sense, time can also be interpreted as the length
parameter of a directed polymer.

\vspace{0.3cm}

A word about random media appellations: we believe the term \textquotedblleft
random environment\textquotedblright\ normally implies that the underlying
randomness is allowed to change over time; the appellation \textquotedblleft
random scenery\textquotedblright\ or \textquotedblleft random
landscape\textquotedblright\ is more specifically used for an environment that
does not change over time; the models we consider herein fall under the
time-varying \textquotedblleft environment\textquotedblright\ umbrella. We now
give some brief specifics about these models.

\vspace{0.3cm}

\noindent\textbf{(1)} We first consider a Brownian polymer in a Gaussian
environment: the polymer itself is modeled by a Brownian motion $b=\{b_{t}%
;t\geq0\}$, defined on a complete filtered probability space $(\mathcal{C}%
,\mathcal{F},(\mathcal{F}_{t})_{t\geq0},(P_{b}^{x})_{x\in{\mathbb{R}}^{d}})$,
where $P_{b}^{x}$ stands for the Wiener measure starting from the initial
condition $x$. The corresponding expected value is denoted by $E_{b}^{x}$, or
simply by $E_{b}$ when $x=0$.

\vspace{0.1cm}

The random environment is represented by a centered Gaussian random field $W$
indexed by ${\mathbb{R}}_{+}\times{\mathbb{R}^{d}}$, defined on another
independent complete probability space $(\Omega,\mathcal{G},\mbox{{\bf P}})$.
Denoting by $\mbox{{\bf E}}$ the expected value with respect to
$\mbox{{\bf P}}$, the covariance structure of $W$ is given by
\begin{equation}
\mbox{{\bf E}}\left[  W(t,x)W(s,y)\right]  =(t\wedge s)\cdot Q(x-y),
\label{E:en2}%
\end{equation}
for a given homogeneous covariance function $Q:{\mathbb{R}}^{d}\rightarrow
{\mathbb{R}}$ satisfying some regularity conditions that will be specified
later on. In particular, the function $t\mapsto\lbrack Q(0)]^{-1/2}W(t,x)$
will be a standard Brownian motion for any fixed $x\in{\mathbb{R}}^{d}$; for
every fixed $t\in{\mathbb{R}}_{+}$, the process $x\mapsto t^{-1/2}W(t,x)$ is a
homogeneous Gaussian field on ${\mathbb{R}}^{d}$ with covariance function $Q$.
Notice that the homogeneity assumption is made here for sake of readability,
but could be weakened for almost all the results we will show. The interested
reader can consult \cite{FV} for the types of tools needed for such generalizations.

\vspace{0.1cm}

Once $b$ and $W$ are defined, the polymer measure itself can be described as
follows: for any $t>0$, the energy of a given path (or configuration) $b$ on
$[0,t]$ is given by the \emph{Hamiltonian}
\begin{equation}
-H_{t}(b)=\int_{0}^{t}W(ds,b_{s}). \label{nrj}%
\end{equation}
A completely rigorous meaning for this integral will be given in the next
section, but for the moment, observe that for any fixed path $b$, $H_{t}(b)$
is a centered Gaussian random variable with variance $tQ(0)$. Based on this
Hamiltonian, for any $x\in{\mathbb{R}}^{d}$, and a given constant $\beta$
(interpreted as the inverse of the temperature of the system), we define our
(random) polymer measure $G_{t}^{x}$ (with $G_{t}:=G_{t}^{0}$) as follows:
\begin{equation}
dG_{t}^{x}(b)=\frac{e^{-\beta H_{t}(b)}}{Z_{t}^{x}}dP_{b}^{x}(b),\quad
\mbox{ with }\quad Z_{t}^{x}=E_{b}^{x}\left[  e^{-\beta H_{t}(b)}\right]  .
\label{defgt}%
\end{equation}

\vspace{0.1cm}

\noindent\textbf{(2)} The second model we consider in this article is the
continuous time random walk on ${\mathbb{Z}}^{d}$ in a white noise potential,
which can be defined similarly to the Brownian polymer above: the polymer is
modeled by a continuous time random walk $\hat{b}=\{\hat{b}_{t};t\geq0\}$ on
${\mathbb{Z}}^{d}$, defined on a complete filtered probability space
$(\hat{\mathcal{C}},\hat{\mathcal{F}},(\hat{\mathcal{F}}_{t})_{t\geq0}%
,(\hat{P}_{\hat{b}}^{x})_{x\in{\mathbb{Z}}^{d}})$. The corresponding expected
value will be denoted by $\hat{E}_{\hat{b}}^{x}$, or simply by $\hat{E}%
_{\hat{b}}$ when $x=0$. Notice that $\hat{b}$ can be represented in terms of
its jump times $\{\tau_{i};i\geq0\}$ and its positions $\{x_{i};i\geq0\}$
between the jumps, as $\hat{b}_{t}=\sum_{i\geq0}x_{i}\mathbf{1}_{[\tau
_{i},\tau_{i+1})}(t)$. Then, under $\hat{P}_{\hat{b}}$, $\tau_{0}=x_{0}=0$,
the sequence $\{\tau_{i+1}-\tau_{i};i\geq0\}$ is i.i.d with common exponential
law $\mathcal{E}(2d)$, and the sequence $\{x_{i};i\geq0\}$ is a nearest
neighbor symmetric random walk on ${\mathbb{Z}}^{d}$.

\vspace{0.1cm}

In this context, the random environment $\hat{W}$ will be defined as a
sequence $\{\hat{W}(.,z);z\in{\mathbb{Z}}^{d}\}$ of Brownian motions, defined
on another independent complete probability space $(\hat{\Omega}%
,\hat{\mathcal{G}},$ $\mathbf{\hat{P}})$. Just like in the Brownian case
described above, the covariance structure we assume on $\hat{W}$ is of the
following type:
\begin{equation}
\mathbf{\hat{E}}\left[  \hat{W}(t,x)\hat{W}(s,y)\right]  =[t\wedge s]\,\hat
{Q}(x-y),\label{E:en3}%
\end{equation}
for a covariance function $\hat{Q}$ defined on ${\mathbb{Z}}^{d}$. Note that
the case where $\hat{Q}\left(  z\right)  =0$ for all $z$ except $\hat
{Q}\left(  0\right)  >0$, is the case where Brownian motions in the family
$\{\hat{W}\left(  \cdot,z\right)  ;z\in{\mathbb{Z}}^{d}\}$ are independent,
i.e. the case of space-time white noise. The Hamiltonian of our system can be
defined formally similarly to the continuous case, as
\[
-\hat{H}_{t}(\hat{b})=\int_{0}^{t}\hat{W}(ds,\hat{b}_{s}).
\]
Notice however that, since $b$ is a piecewise constant function, the
Hamiltonian $\hat{H}_{t}(\hat{b})$ can also be written as
\begin{equation}
-\hat{H}_{t}(\hat{b})=\sum_{i=0}^{N_{t}}\hat{W}(\tau_{i+1},x_{i})-\hat{W}%
(\tau_{i},x_{i}),\label{eq:expression-hamil-discrete}%
\end{equation}
where $N_{t}$ designates the number of jumps of $\hat{b}$ before time $t$, and
$\tau_{N_{t}+1}=t$ by convention. Once the Hamiltonian $\hat{H}_{t}$ is
defined, a Gibbs-type measure $\hat{G}_{t}$ can be introduced similarly to
(\ref{defgt}) in the Brownian case.\vspace{0.3cm}

As mentioned before, our aim in this article is to give some sharp estimates
on the free energies $p(\beta)$ and $\hat{p}(\beta)$ of the two systems
described above, for large $\beta$. The quantities of interest are defined
asymptotically as
\[
p(\beta)=\lim_{t\rightarrow\infty}\frac{1}{t}\mbox{{\bf E}}\left[  \log
(Z_{t})\right]  ,\quad\mbox{ and }\quad\hat{p}(\beta)=\lim_{t\rightarrow
\infty}\frac{1}{t}\mbox{{\bf E}}\left[  \log(\hat{Z}_{t})\right]  ;
\]
it is well-known (see e.g. \cite{RT} for the Brownian case) that the limits
above exist, are typically positive, and are both bounded from above by
$Q(0)\beta^{2}/2$. It is then possible to separate a region of weak disorder
from a region of strong disorder according to the value of $p(\beta)$: we will
say that the polymer is in the weak disorder regime if $p(\beta)=Q(0)\beta
^{2}/2$, while the strong disorder regime is defined by the strict inequality
$p(\beta)<Q(0)\beta^{2}/2$. These two notions have some nice interpretations
in terms of the behavior of the particle under the Gibbs measure (see e.g.
\cite{CH2,CY3}), and it is expected, for any model of polymer in a random
environment, that the strong disorder regime is attained whenever $\beta$ is
large enough. It is then natural to ask if one can obtain a sharper
information than $p(\beta)<Q(0)\beta^{2}/2$ in the low temperature phase.
Indeed, on the one hand, this may quantify in a sense how far we are from the
weak disorder regime, and how much localization there is on our measures
$G_{t},\hat{G}_{t}$. On the other hand, the penalization method explained in
\cite{RY} can be roughly summarized in the following way: if one can get a
sharp equivalent for the quantity $E_{b}[e^{-\beta H_{t}(b)}]$, then this will
also allow a detailed description of the limit $\lim_{t\rightarrow\infty}%
G_{t}$. This latter program is of course beyond the scope of the current
article, but is a good motivation for getting some precise information about
the function $p(\beta)$.

\subsection{Summary of results}

We now describe our main results. Our principal result in continuous space
will be obtained in terms of the regularity of $Q$ in a neighborhood of 0. In
particular, we shall assume some upper and lower bounds on $Q$ of the form
\begin{equation}
c_{0}|x|^{2H}\leq Q(0)-Q(x)\leq c_{1}|x|^{2H},\quad
\mbox{ for all }x\mbox{ such that }|x|\in\lbrack0,r_{0}], \label{Holderhypo}%
\end{equation}
for a given exponent $H\in(0,1]$ and $r_{0}>0$. It should be noticed that
condition (\ref{Holderhypo}) is equivalent to assuming that $W$ has a specific
almost-sure modulus of continuity in space, of order $\left\vert x\right\vert
^{H}\log^{1/2}\left(  1/\left\vert x\right\vert \right)  $, i.e. barely
failing to be $H$-H\"{o}lder continuous (see \cite{TTV} for details). Then,
under these conditions, we will get the following conclusions.

\begin{theorem}
\label{thm:bnd-free-nrj-brown-holder} Assume that the function $Q$ satisfies
condition (\ref{Holderhypo}). Then the following hold true:

\begin{enumerate}
\item If $H\in\lbrack1/2,1]$, we have for some constants $C_{0,d}$ and
$C_{1,d}$ depending only on $Q$ and $d$, for all $\beta\geq1$,
\[
C_{0,d}\beta^{4/3}\leq p(\beta) \leq C_{1,d}\beta^{2-2H/\left(  3H+1\right)
}.
\]

\item If $H\in(0,1/2]$, we have for some constants $\beta_{Q}$, $C_{0,d}%
^{\prime}$, and $C_{1,d}^{\prime}$ depending only on $Q$ and $d$, for all
$\beta\geq\beta_{Q}$,
\[
C_{0,d}^{\prime}\beta^{2/(1+H)}\leq p(\beta)\leq C_{1,d}^{\prime}%
\beta^{2-2H/\left(  3H+1\right)  }.
\]

\end{enumerate}
\end{theorem}

\begin{figure}[ht]
\label{fig:exp-beta}
\par
\begin{center}
\includegraphics[scale=1.0]{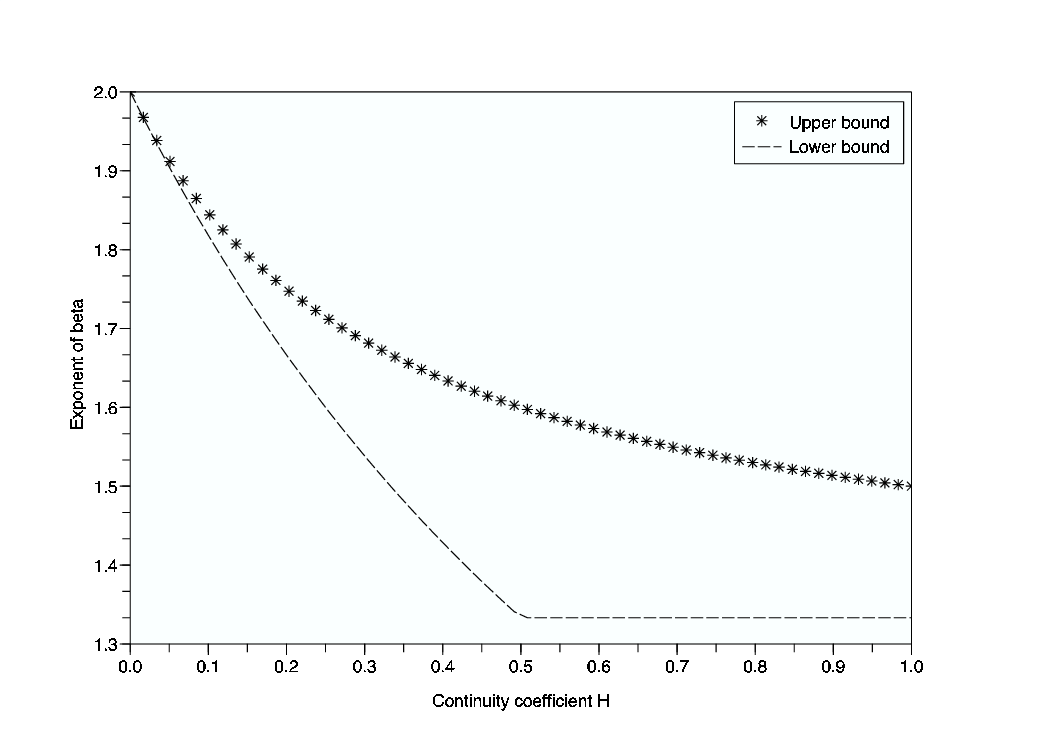}
\end{center}
\caption{Exponent of $\beta$ in function of $H$}%
\end{figure}

Corresponding almost sure results on $t^{-1}\mbox{{\bf E}}\left[  \log
(Z_{t})\right]  $ also hold, as seen in Corollary \ref{cor2} and Proposition
\ref{p14} below. Let us make a few elementary comments about the above
theorem's bounds, which are also summarized in Figure \ref{fig:exp-beta}.
First of all, the exponent of $\beta$ in those estimates is decreasing with
$H$, which seems to indicate a stronger disorder when the Gaussian field $W$
is smoother in space. Furthermore, in the case $H\in\lbrack1/2,1]$, the gap
between the two estimates decreases as $H$ increases to $1$; for $H=1/2$, we
get bounds with the powers of $\beta$ equal to $4/3$ and $8/5;$ and for $H=1$,
the bounds are $4/3$ and $3/2$. It should be noted that the case $H=1/2$ is
our least sharp result, while the case $H=1$ yields the lowest power of
$\beta$; one should not expect lower powers for any potential $W$ even if $W$
is so smooth that it is $C^{\infty}$ in space: indeed, unless $W$ is highly
degenerate, the lower bound in (\ref{Holderhypo}) should hold with $H=1$,
while the upper bound will automatically be satisfied with $H=1$. The case of
small $H$ is more interesting. Indeed, we can rewrite the lower and upper
bounds above as
\[
C_{0,d}^{\prime}\beta^{2-2H+F\left(  H\right)  }\leq p(\beta)\leq
C_{1,d}^{\prime}\beta^{2-2H+G\left(  H\right)  }%
\]
where the functions $F$ and $G$ satisfy $F\left(  x\right)  =2x^{2}+O\left(
x^{3}\right)  \ $and $G\left(  x\right)  =6x^{2}+O\left(  x^{3}\right)  $ for
$x$ near $0$. We therefore see that the asymptotic $\beta^{2-2H}$ is quite
sharp for small $H$, but that the second order term in the expansion of the
power of $\beta$ for small $H$, while bounded, is always positive.

Using ideas introduced in \cite{FV} to deal with spatially non-homogeneous
media, it is possible to extend Theorem \ref{thm:bnd-free-nrj-brown-holder}.
The reader will check that the first of the following two corollaries is
trivial to prove using the tools in this article. The second corollary
requires techniques in \cite{FV}, and can also be proved directly by using
sub-Gaussian concentration results (see \cite{VV}). We do not give any details
of its proof, for the sake of conciseness. Neither corollary assumes that $W$
is spatially homogeneous. One will note that no assertion on the existence of
$p\left(  \beta\right)  $ is made in these corollaries, but that the first
corollary already implies strong disorder for large $\beta$ in the sense that
$\limsup t^{-1}\mbox{{\bf E}}\left[  \log(Z_{t})\right]  <\beta^{2}Q\left(
0\right)  /2$. \cite{FV} can be consulted for conditions under which $p\left(
\beta\right)  $ exists even if $W$ is not spatially homogeneous.

\begin{corollary}
\label{cor1} In the non homogeneous case, the following bounds are satisfied:

\begin{itemize}
\item \emph{[Upper bound] }Assume that for some $r_{0},c_{1}>0$, for all
$x,y\in{\mathbb{R}^{d}}$ such that $\left\vert x-y\right\vert \leq r_{1}$, the
spatial canonical metric of $W$ is bounded above as%
\[
\delta^{2}\left(  x,y\right)  :=\mathbf{E}\left[  \left(  W\left(  1,x\right)
-W\left(  1,y\right)  \right)  ^{2}\right]  \leq c_{1}\left\vert
x-y\right\vert ^{2H}.
\]
Then, replacing $p\left(  \beta\right)  $ by $\limsup_{\beta\rightarrow\infty
}t^{-1}\mbox{{\bf E}}\left[  \log(Z_{t})\right]  $, the two upper bound
results in Theorem \ref{thm:bnd-free-nrj-brown-holder} hold.

\item \emph{[Lower bound] }Assume that for some $r_{0},c_{0}>0$, for all
$x,y\in{\mathbb{R}^{d}}$ such that $\left\vert x-y\right\vert \leq r_{0}$, we
have%
\[
\delta^{2}\left(  x,y\right)  :=\mathbf{E}\left[  \left(  W\left(  1,x\right)
-W\left(  1,y\right)  \right)  ^{2}\right]  \geq c_{0}\left\vert
x-y\right\vert ^{2H}.
\]
Then, replacing $p\left(  \beta\right)  $ by $\liminf_{\beta\rightarrow\infty
}t^{-1}\mbox{{\bf E}}\left[  \log(Z_{t})\right]  $, the two lower bound
results in Theorem \ref{thm:bnd-free-nrj-brown-holder} hold.
\end{itemize}
\end{corollary}

\begin{corollary}
\label{cor2} Under the hypotheses of Corollary \ref{cor1}, its conclusions
also hold $\mathbf{P}$-almost surely with $\limsup_{\beta\rightarrow\infty
}t^{-1}\mbox{{\bf E}}\left[  \log(Z_{t})\right]  $ replaced by $\limsup
_{\beta\rightarrow\infty}t^{-1}\log(Z_{t})$, and similarly for the $\liminf$ 's.
\end{corollary}

Since our estimates become sharper as $H\rightarrow0$, and also due to the
fact that the behavior of $p(\beta)$ is nearly quadratic in $\beta$ for small
$H$ (i.e. approaching the weak disorder regime), we decided to explore further
the region of logarithmic spatial regularity for $W$, in order to determine
whether one ever leaves the strong disorder regime. Namely, we also examine
the situation of a covariance function $Q$ for which there exist positive
constants $c_{0}$, $c_{1}$, and $r_{1}$ such that for all $x$ with $\left\vert
x\right\vert \leq r_{1}$,%
\begin{equation}
c_{0}\log^{-2\gamma}\left(  1/\left\vert x\right\vert \right)  \leq
Q(0)-Q(x)\leq c_{1}\log^{-2\gamma}\left(  1/\left\vert x\right\vert \right)  ,
\label{loghypo}%
\end{equation}
where $\gamma$ is a given positive exponent. Assumption (\ref{loghypo})
implies that $W$ is not spatially H\"{o}lder-continuous for any exponent
$H\in(0,1]$. Moreover, the theory of Gaussian regularity implies that, if
$\gamma>1/2$, $W$ is almost-surely continuous in space, with modulus of
continuity proportional to $\log^{-\gamma+1/2}\left(  1/\left\vert
x\right\vert \right)  $, while if $\gamma\leq1/2$, $W$ is almost-surely not
uniformly continuous on any interval in space, and in fact is unbounded on any
interval. We will then establish the following result, which is optimal, up to
multiplicative constants.

\begin{theorem}
\label{thm:bnd-free-nrj-brown-log} Assume condition (\ref{loghypo}) where
$\gamma>0$. We have for some constants $D_{0,d}$ and $D_{1,d}$ depending only
on $Q$ and $d$, for all $\beta$ large enough,
\[
D_{0,d}\beta^{2}\log^{-2\gamma}\left(  \beta\right)  \leq p(\beta) \leq
D_{1,d}\beta^{2}\log^{-2\gamma}\left(  \beta\right)  .
\]

\end{theorem}

Besides giving a sharp result up to constants for the free energy $p(\beta)$,
the last result will allow us to make a link between our Brownian model and
the random walk polymer described by the Hamiltonian
(\ref{eq:expression-hamil-discrete}). Indeed, the following result will also
be proved in the sequel.

\begin{theorem}
\label{thm:bnd-free-nrj-walk} Assume that $\hat Q(0)-\hat Q(2)>0$, where $\hat
Q$ has been defined at (\ref{E:en3}). Then the free energy $\hat p(\beta)$ of
the random walk polymer $\hat b$ satisfies, for $\beta$ large enough:
\begin{equation}
\label{eq:bound-p-beta-discrete}D^{\prime}_{0,d}\beta^{2}\log^{-1}\left(
\beta\right)  \leq\hat p(\beta) \leq D^{\prime}_{1,d}\beta^{2}\log^{-1}\left(
\beta\right)  ,
\end{equation}
for two constants $D_{0,d}$ and $D_{1,d}$ depending only on $Q$ and $d$.
\end{theorem}

Relation (\ref{eq:bound-p-beta-discrete}) will be obtained here thanks to some
simple arguments, which allow the extension to spatially inhomogeneous media.
In the special homogeneous case of space-time white noise ($Q\left(  x\right)
=0$ for all $x\neq0$), more can be said: the exact value of the limit
$\lim_{\beta\rightarrow\infty}\hat{p}(\beta)\log\left(  \beta\right)
/\beta^{2}$ can be computed in this situation; this result has been
established by the authors of the work in preparation \cite{CMRT}.

In relation with the continuous space model considered at Theorem
\ref{thm:bnd-free-nrj-brown-log}, we see that to obtain the same behavior as
with space-time white noise in discrete space, we need to use precisely the
environment $W$ in ${\mathbb{R}^{d}}$ with the logarithmic regularity
corresponding to $\gamma=1/2$ in (\ref{loghypo}). As mentioned before, this
behavior of $W$ happens to be exactly at the threshold in which $W$ becomes
almost-surely discontinuous and unbounded on every interval. Nevertheless such
a $W$ is still function-valued. Hence, for the purpose of understanding the
polymer partition function, there is no need to study the space-time white
noise in continuous space, for which $W\left(  t,\cdot\right)  $ is not a
bonafide function (only a distribution), and for which the meaning of $Z_{t}$
itself is difficult to even define. Another way to interpret the coincidence
of behaviors for \textquotedblleft\emph{space-time white noise in
}${\mathbb{R}}_{+}\times{\mathbb{Z}}^{d}$\textquotedblright\ and for
\textquotedblleft$\gamma=1/2$\textquotedblright\ is to say that both models
for $W$ are function-valued and exhibit spatial discontinuity: indeed, in
discrete space, one extends $W\left(  t,\cdot\right)  $ to ${\mathbb{R}}^{d}$
by making it piecewise constant, in order to preserve independence. The fact
that the limit in Theorem \ref{thm:bnd-free-nrj-brown-log} depends on $\gamma$
does prove, however, that the continuous-space polymer model under logarithmic
regularity is richer than the discrete-space one.

As in the H\"older-scale continuous space setting, we have the following
corollaries, in which $W$ is allowed to be spatially inhomogeneous. Again, we
do not include proofs of these results for the sake of conciseness.

\begin{corollary}
Assume the lower and upper bound hypotheses in Corollary \ref{cor1} hold with
$\left\vert x-y\right\vert ^{2H}$ replaced by $\log^{-2\gamma}\left(
1/\left\vert x-y\right\vert \right)  $. Then the conclusions of Theorem
\ref{thm:bnd-free-nrj-brown-log} hold with $p\left(  \beta\right)  $ replaced
by $\liminf_{\beta\rightarrow\infty}t^{-1}\mbox{{\bf E}}\left[  \log
(Z_{t})\right]  $ for the lower bound, and by $\limsup_{\beta\rightarrow
\infty}t^{-1}\mbox{{\bf E}}\left[  \log(Z_{t})\right]  $ for the upper bound.
Almost-sure results as in Corollary \ref{cor2} also hold.
\end{corollary}

\begin{corollary}
For the discrete space polymer in Theorem \ref{thm:bnd-free-nrj-walk}, assume,
instead of $\hat{Q}\left(  0\right)  >\hat{Q}\left(  2\right)  $, that
$\mathbf{E}\left[  \left(  W\left(  1,-1\right)  -W\left(  1,1\right)
\right)  ^{2}\right]  >0$. Then the conclusions of Theorem
\ref{thm:bnd-free-nrj-walk} hold with $\hat{p}\left(  \beta\right)  $ replaced
by $\liminf_{\beta\rightarrow\infty}t^{-1}\mathbf{\hat{E}}[\log(\hat{Z}_{t})]$
for the lower bound, and with $\hat{p}\left(  \beta\right)  $ replaced by
$\limsup_{\beta\rightarrow\infty}t^{-1}\mathbf{\hat{E}}[\log(\hat{Z}_{t})]$
for the upper bound. Almost-sure results as in Corollary \ref{cor2} also
hold.\bigskip
\end{corollary}

Let us say a few words now about the methodology we have used in order to get
our results. It is inspired by the literature on Lyapounov exponents for
stochastic {\small PDE}s \cite{CKM,CM,CV,FV,TV1,TV2}; our upper bound results
rely heavily on the estimation of the supremum of some well-chosen Gaussian
fields, using such results as Dudley's so-called entropy upper bound, and the
Borell-Sudakov inequality (see \cite{A} or \cite{VV}); our lower bound results
are obtained more \textquotedblleft by hand\textquotedblright, by isolating
very simple polymer configurations $b$ or $\hat{b}$ which maximize the random
medium's increments in the Hamiltonian $H_{t}\left(  b\right)  $ or $\hat
{H}_{t}( \hat{b}) $, and showing that these configurations contain enough
weight to provide lower bounds. It turns out that these estimation procedures
works better when the configuration $b$ is simple enough, such as a piecewise
constant or linear function. For the upper bound in the continuous case, a
careful discretization of our Brownian path will thus have to be performed in
order to get our main results; the resulting proof cannot exploit the discrete
case itself because of the different nature of the discrete and continuous environments.

\vspace{0.3cm}

The structure of the article is as follows: Section 2 contains preliminary
information on the partition function. Section 3 deals with the Brownian
polymer. Section 4 covers the random walk polymer. In order to simplify the
notation, throughout the paper we will use $C$ to represent the constants, but
acknowledge that the value it represents will change, even from line to line.

\section{Preliminaries; the partition function}

In this section, we will first recall some basic facts about the definition
and the simplest properties of the partition functions $Z_{t}$ and $\hat
{Z}_{t}$ which have been already considered in the introduction. We will also
give briefly some notions of Gaussian analysis which will be used later
on.\vspace{0.3cm}

We begin with basic information about the partition function of the Brownian
polymer. Recall that $W$ is a centered Gaussian field on ${\mathbb{R}}%
_{+}\times{\mathbb{R}}^{d}$, defined by its covariance structure
(\ref{E:en2}). The Hamiltonian $H_{t}(b)$ given by (\ref{nrj}) can be defined
more rigorously through a Fourier transform procedure: there exists (see e.g.
\cite{CV} for further details) a centered Gaussian independently scattered
$\mathbb{C}$-valued measure $\nu$ on ${\mathbb{R}}_{+}\times{\mathbb{R}}^{d}$
such that
\begin{equation}
W(t,x)=\int_{{\mathbb{R}}_{+}\times{\mathbb{R}}^{d}}\mathbf{1}_{[0,t]}%
(s)e^{\imath ux}\nu(ds,du), \label{E:rep_int}%
\end{equation}
where the simple notation $ux$ stands for the inner product $u\cdot x$ in
${\mathbb{R}}^{d}$. For every test function $f:{\mathbb{R}}_{+}\times
{\mathbb{R}}^{d}\rightarrow\mathbb{C}$, set now
\begin{equation}
\nu(f)\equiv\int_{{\mathbb{R}}_{+}\times{\mathbb{R}}^{d}}f(s,u)\nu(ds,du).
\label{defnuf}%
\end{equation}
While the random variable $\nu\left(  f\right)  $ may be complex-valued, to
ensure that it is real valued, it is sufficient to assume that $f$ is of the
form $f\left(  s,u\right)  =f_{1}\left(  s\right)  e^{\imath uf_{2}\left(
s\right)  }$ for real valued functions $f_{1}$ and $f_{2}$. Then the law of
$\nu$ is defined by the following covariance structure: for any such test
functions $f,g:{\mathbb{R}}_{+}\times{\mathbb{R}}^{d}\rightarrow\mathbb{C}$,
we have
\begin{equation}
\mbox{{\bf E}}\left[  \nu(f)\nu(g)\right]  =\int_{{\mathbb{R}}_{+}%
\times{\mathbb{R}}^{d}}f(s,u)\overline{g(s,u)}\hat{Q}(du)ds, \label{E:cov_X}%
\end{equation}
where the finite positive measure $\hat{Q}$ is the Fourier transform of $Q$
(see \cite{tv99} for details).

From (\ref{E:rep_int}), we see that the It\^{o}-stochastic differential of $W$
in time can be understood as $W\left(  ds,x\right)  :=\int_{u\in{\mathbb{R}%
}^{d}}e^{\imath ux}\nu(ds,du)$, or even, if the measure $\hat{Q}\left(
du\right)  $ has a density $f\left(  u\right)  $ with respect to the Lebesgue
measure, which is typical, as
\[
W\left(  ds,x\right)  :=\int_{u\in{\mathbb{R}}^{d}}e^{\imath ux}\sqrt{f\left(
u\right)  }M(ds,du)
\]
where $M$ is a white-noise measure on ${\mathbb{R}}_{+}\times{\mathbb{R}}^{d}%
$, i.e. a centered independently scattered Gaussian measure with covariance
given by $\mathbf{E}\left[  M\left(  A\right)  M\left(  B\right)  \right]
=m_{Leb}\left(  A\cap B\right)  $ where $m_{Leb}$ is Lebesgue's measure on
${\mathbb{R}}_{+}\times{\mathbb{R}}^{d}$.

We can go back now to the definition of $H_{t}(b)$: invoking the
representation (\ref{E:rep_int}), we can write
\begin{equation}
-H_{t}(b):=\int_{0}^{t}W(ds,b_{s})=\int_{0}^{t}\int_{{\mathbb{R}}^{d}%
}e^{\imath ub_{s}}\nu(ds,du),\label{E:rep_intH}%
\end{equation}
taking this expression as a definition of $H_{t}\left(  b\right)  $ for each
fixed path $b$; it can be shown (see \cite{CV}) that the right hand side of
the above relation is well defined for any H\"{o}lder continuous path $b$, by
a $L^{2}$-limit procedure. Such a limiting procedure can be adapted to the
specific case of constructing $H_{t}\left(  b\right)  $, using the natural
time evolution structure; we will not comment on this further. However, the
reader will surmise that the following remark, given for the sake of
illustration, can be useful: when $\hat{Q}$ has a density $f$, we obtain
$-H_{t}(b)=\iint_{[0,t]\times{\mathbb{R}}^{d}}e^{\imath ub_{s}}\sqrt{f\left(
u\right)  }M\left(  ds,du\right)  .\vspace*{0.1in}$

We use as the definition of the partition function $Z_{t}^{x}$, its expression
in (\ref{defgt}), and set its expectation under $\mathbf{P}$ as%
\begin{equation}
p_{t}(\beta):=\frac{1}{t}\mbox{{\bf E}}\left[  \log\left(  Z_{t}^{x}\right)
\right]  , \label{fren}%
\end{equation}
usually called the free energy of the system. It is easily seen that
$p_{t}(\beta)$ is independent of the initial condition $x\in{\mathbb{R}}^{d}$,
thanks to the spatial homogeneity of $W$. Thus, in the remainder of the paper,
$x$ will be understood as $0$ when not specified, and $E_{b},Z_{t}$ will stand
for $E_{b}^{0},Z_{t}^{0}$, etc... We summarize some basic results on
$p_{t}(\beta)$ and $Z_{t}$ established in \cite{RT}.

\begin{proposition}
\label{p14} For all $\beta>0$ there exists a constant $p(\beta)>0$ such that
\begin{equation}
\label{eq:lim-p-t-beta}p(\beta):=\lim_{t\rightarrow\infty}p_{t}(\beta
)=\sup_{t\geq0}p_{t}(\beta).
\end{equation}
Furthermore, the function $p$ satisfies:

\begin{enumerate}
\item The map $\beta\mapsto p(\beta)$ is a convex nondecreasing function on
${\mathbb{R}}_{+}$.

\item The following upper bound holds true:
\begin{equation}
\label{roughbnd}p(\beta)\le\frac{\beta^{2}}{2} Q(0).
\end{equation}

\item $\mbox{{\bf P}}$-almost surely, we have
\begin{equation}
\lim_{t\rightarrow\infty}\frac{1}{t}\log Z_{t}=p(\beta).\label{aslimfree}%
\end{equation}

\end{enumerate}
\end{proposition}

For the random walk polymer on ${\mathbb{Z}}^{d}$, the Hamiltonian $\hat
{H}_{t}(\hat{b})$ is easier to define, and can be expressed in a simple way by
(\ref{eq:expression-hamil-discrete}). Recall then that $\hat{Z}_{t},\hat
{p}_{t}(\beta)$ are defined as:
\[
\hat{Z}_{t}=\hat{E}_{\hat{b}}\left[  e^{-\hat{H}_{t}(\hat{b})}\right]
,\quad\mbox{ and }\quad\hat{p}_{t}(\beta)=\hat{\mbox{{\bf E}}}\left[
\log(\hat{Z}_{t})\right]  .
\]
Then, using the same kind of arguments as in \cite{RT} (see also \cite{CMRT}),
we get the following:

\begin{proposition}
\label{prop:simple-prop-free-discrete} The same conclusions as in Proposition
\ref{p14} hold true for the random walk polymer $\hat b$.
\end{proposition}

\section{Estimates of the free energy: continuous space}

In this section, we will proceed to the proof of Theorems
\ref{thm:bnd-free-nrj-brown-holder} and, \ref{thm:bnd-free-nrj-brown-log}, by
means of some estimates for some well-chosen Gaussian random fields.

The hypothesis we use guarantees that there is some $H\in(0,1)$ such that $W$
is no more than $H$-H\"{o}lder continuous in space. Accordingly, we define the
homogeneous spatial \textit{canonical metric} $\delta$ of $W$ by
\begin{equation}
\delta^{2}(x-y):=\mathbf{E}\left[  \left(  W(1,x)-W(1,y)\right)  ^{2}\right]
=2\left(  Q(0)-Q(x-y)\right)  , \label{eq:def-delta}%
\end{equation}
for all $x,y\in\mathbb{R}^{d}$. Our hypotheses on $\delta$ translate
immediately into statements about $Q$ via this formula.

In our results below, we have also tried to specify the dependence of our
constants on the dimension of the space variable. An interesting point in that
respect is given in the lower bound of Subsection \ref{sec:low-bnd-brownian}
below, which has to do with weak versus strong disorder in very
high-dimensional cases.

\subsection{Upper bound in the Brownian case}

The upper bound in Theorem \ref{thm:bnd-free-nrj-brown-holder} follows
immediately from the following proposition, which proves in particular that
strong disorder holds for all $H\in(0,1]$.

\begin{proposition}
\label{prop:upper-bnd-holder} Assume that there exist a number $H\in(0,1]$ and
numbers $c_{1},r_{1}$ such that for all $x,y\in\mathbb{R}^{d}$ with
$\left\vert x-y\right\vert \leq r_{1}$ we have
\begin{equation}
\delta(x-y)<c_{1}\left\vert x-y\right\vert ^{H}. \label{hub}%
\end{equation}
Then there exists a constant $C$ depending only on $Q$ and a constant
$\beta_{0}$ depending only on $r_{1}$ and $d$, such that for all $\beta
\geq\beta_{0}$,%
\[
p(\beta)\leq Cd^{\frac{7H}{1+3H}}\beta^{\frac{2+4H}{1+3H}}.
\]

\end{proposition}

\begin{proof}
Let us divide the proof in several steps:

\vspace{0.3cm}

\noindent\textit{Step 1: Strategy.} From relation (\ref{eq:lim-p-t-beta}), we
have
\[
p(\beta)\leq\limsup_{t\rightarrow\infty}p_{t}(\beta).
\]
Our strategy is then to give an estimation of $p_{t}(\beta)$ for a
\textsl{discretized} path $\tilde{b}\in\varepsilon\mathbb{Z}^{d}$ that stays
close to $b$ and proceeds only by jumps. Thanks to this substitution, and
using H\"{o}lder's and Jensen's inequalities, we shall obtain
\begin{align}
\mathbf{E}\left[  \log(Z_{t})\right]   &  =\mathbf{E}\left[  \log E_{b}\left[
\exp\left(  -\beta\left[  H_{t}(b)-H_{t}(\tilde{b})\right]  \right)
\exp-\beta H_{t}(\tilde{b})\right]  \right]
\label{eq:up-bnd-holder-with-error}\\
&  \leq\frac{1}{2}\mathbf{E}\left[  \log E_{b}\left[  \exp\left(
-2\beta\lbrack H_{t}(b)-H_{t}(\tilde{b})]\right)  \right]  \right]  +\frac
{1}{2}\mathbf{E}\left[  \log E_{b}\left[  \exp\left(  -2\beta H_{t}(\tilde
{b})\right)  \right]  \right]  \nonumber\\
&  \leq\frac{1}{2}\log E_{b}\left[  \exp2\beta^{2}\int_{0}^{t}\left(
\delta(b_{s}-\tilde{b_{s}})\right)  ^{2}ds\right]  +\frac{1}{2}\mathbf{E}%
\left[  \log E_{b}\left[  \exp\left(  -2\beta H_{t}(\tilde{b})\right)
\right]  \right]  .\nonumber
\end{align}
Notice that the first term on the right-hand side represents the error made by
considering the discretized path $\tilde{b}$ instead of $b$, but thanks to
hypothesis (\ref{hub}) and the definition of $\tilde{b}$ we will easily
control it.\vspace{0.3cm}

\noindent\textit{Step 2: The discretized path.} Let us describe now the
discretized process we shall use in the sequel: we will approximate the
Brownian path $b$ with a path that stays in $\varepsilon\mathbb{Z}^{d}$, where
$\varepsilon$ is a small positive number. Let $b^{j}$ be the $j$-th component
of the $d$-dimensional path $b$. Let $T^{j}_{1}$ be the first time that
$b^{j}$ exits the interval $(-\varepsilon,\varepsilon)$ and $T_{i+1}^{j}$ be
the first time after $T_{i}^{j}$ that $b^{j}$ exits $(b_{T_{i}^{j}%
}-\varepsilon,b_{T_{i}^{j}}+\varepsilon)$. So, for a fixed component $j$, the
times $(T_{i+1}^{j}-T_{i}^{j})_{i=0}^{\infty}$ are i.i.d. and the successive
positions $x_{m}^{j}=b_{T_{m}^{j}}^{j}$, which are independent of the jump
times, form a one-dimensional symmetric random walk on $\varepsilon\mathbb{Z}$
in discrete time.

Now let $(T_{n})_{n=0}^{\infty}$ be the increasing sequence of all the
$(T_{m}^{j})_{j,m}$ and let $(x_{n})_{n=0}^{\infty}$ be the nearest neighbor
path in $\varepsilon\mathbb{Z}^{d}$ with $x_{0}=0$ whose $j$-th component
takes the same step as $x_{m}^{j}$ at time $T_{m}^{j}$. We define the
\textsl{discretized path} $\tilde{b}$ as the path that jumps to site $x_{n}$
at time $T_{n}$ and it is constant between jumps.

\begin{remark}
\label{sost} At any time $s$, each coordinate of $\tilde{b}_{s}$ is within
$\varepsilon$ of the corresponding one of $b_{s}$. So the distance separating
the two paths is never more than $\varepsilon\sqrt{d}$. Thus we have, for all
$s\geq0$, $\vert b_{s}-\tilde{b}_{s}\vert\leq\varepsilon d^{1/2}.$
\end{remark}

\begin{remark}
\label{rmk:bdn-error-pt-pt-ep} Thanks to Remark \ref{sost} we can now control
the error term we have defined at relation (\ref{eq:up-bnd-holder-with-error}%
). In fact, owing to Hypothesis (\ref{hub}), we have
\begin{multline*}
\frac{1}{2t}\log E_{b}\left[  \exp2\beta^{2}\int_{0}^{t}\left(  \delta
(b_{s}-\tilde{b_{s}})\right)  ^{2}ds\right] \\
\leq\frac{1}{2t}\log E_{b}\left[  \exp\left(  2\beta^{2}C^{2}\int_{0}%
^{t}|b_{s}-\tilde{b_{s}}|^{2H}dt\right)  \right]  \leq C\beta^{2}%
\varepsilon^{2H}d^{H},
\end{multline*}
where we recall that $C$ is a constant depending on $Q$ that can change from
line to line.
\end{remark}

Plugging this last inequality into (\ref{eq:up-bnd-holder-with-error}), and
defining
\[
p_{t}^{\varepsilon}(\beta)=\frac{1}{t}\mathbf{E}\left[  \log E_{b}\left[
\exp\left(  -2\beta H_{t}(\tilde{b})\right)  \right]  \right]  ,
\]
we have thus obtained the following estimate for $p_{t}(\beta)$:
\begin{equation}
p_{t}(\beta)\leq C\beta^{2}\varepsilon^{2H}d^{H}+\frac{1}{2}p_{t}%
^{\varepsilon}(\beta).\label{pt}%
\end{equation}
We shall try now to get some suitable bounds on $p_{t}^{\varepsilon}(\beta
)$.\vspace{0.3cm}

\noindent\textit{Step 3: Study of }$p_{t}^{\varepsilon}(\beta)$\textit{.} Let
$N_{t}^{j}$ be the number of jumps of the $j$-th component of $\tilde{b}$ up
to time $t$. For a multi-index $k=(k_{1},\cdots,k_{d})$ let $|k|=k_{1}%
+\cdots+k_{d}$, so the total number of jumps of $\tilde{b}$ up to time $t$ is
$|N_{t}|=N_{t}^{1}+\cdots+N_{t}^{d}$. Denote by $\mathcal{S}(t,n)$ the simplex
of all possible sequences of $n$ jump times up to time $t$, namely
\begin{equation}
\mathcal{S}(t,n)=\left\{  \mathbf{t}=\left(  t_{0},\cdots,t_{n}\right)
:0=t_{0}\leq\cdots\leq t_{n}\leq t\right\}  . \label{eq:def-simplex}%
\end{equation}
The set of the first $k_{j}$ jump times of the $j$-th component of $\tilde{b}$
is a point $(t_{i}^{j})_{i=1}^{k_{j}}$ in $\mathcal{S}(t,k_{j})$. Given the
set of all jump times $\left\{  t_{i}^{j}:j\in\left[  1,\cdots,d\right]
;i\in\left[  1,\cdots,k_{j}\right]  \right\}  $, let $\left\{  \tilde{t}%
_{l}:l\in\left[  0,|k|+1\right]  \right\}  $ be the same set but ordered and
with the convention $\tilde{t}_{0}=0$, $\tilde{t}_{|k|+1}=t$. And finally let
$\tilde{x}_{l}$ be the value of $\tilde{b}$ between the two jump times
$\tilde{t}_{l}$ and $\tilde{t}_{l+1}$. Denote by $\mathcal{P}_{n}$ the set of
all such $\tilde{x}$, i.e. the set of all nearest-neighbor random walk paths
of length $k$ starting at the origin.

Then if we fix $|N_{t}|=|k|$, we can write
\[
H_{t}(\tilde{b})=X\left(  |k|,(\tilde{t}_{l})_{l=1}^{|k|},(\tilde{x}%
_{l})_{l=1}^{|k|}\right)  ,
\]
where
\[
X\left(  |k|,(\tilde{t}_{l})_{l=1}^{|k|},(\tilde{x}_{l})_{l=1}^{|k|}\right)
=\sum_{i=0}^{|k|}\left[  W(\tilde{t}_{i+1},\tilde{x}_{i})-W(\tilde{t}%
_{i},\tilde{x}_{i})\right]  .
\]
Thanks to these notations, we have
\begin{align*}
tp_{t}^{\varepsilon}(\beta)  &  =\mathbf{E}\left[  \log E_{b}\left[
\exp(-2\beta H_{t}(\tilde{b}))\right]  \right] \\
&  =\mathbf{E}\left[  \log E_{b}\left[  \exp\left(  -2\beta X\left(
|N_{t}|,(\tilde{t}_{l})_{l=1}^{|N_{t}|},(\tilde{x}_{l})_{l=1}^{|N_{t}%
|}\right)  \right)  \right]  \right]  .
\end{align*}
So we can write the expectation with respect to $b$ as:
\begin{align*}
E_{b}\left[  \exp(-2\beta H_{t}(\tilde{b}))\right]   &  =\sum_{n\geq1}%
E_{b}\left[  \exp(-2\beta H_{t}(\tilde{b}))\Big|\Big.|N_{t}|\in\left[
t\alpha(n-1),t\alpha n\right]  \right] \\
&  \hspace{2cm}\times P_{b}\left[  |N_{t}|\in\left[  t\alpha(n-1),t\alpha
n\right]  \right]  .
\end{align*}
The number of jumps of the discretized path $\tilde{b}$ in a given interval
$[0,t]$ will play a crucial role in our optimization procedure. For a
parameter $\alpha>0$ which will be fixed later on, let us thus define
\[
T_{n\alpha}=\left\{  \left(  k,\tilde{t},\tilde{x}\right)  :k\leq
tn\alpha,\tilde{t}\in\mathcal{S}(t,k),\tilde{x}\in\mathcal{P}_{k}\right\}  .
\]
Then the following estimates will be essential for our future computations:
\begin{align}
P_{b}\left[  N_{t}^{j}>n\alpha t\right]   &  \leq\exp\left(  -\frac{t}%
{2}(\alpha n\varepsilon)^{2}+t\alpha n\right) \label{stimaprob}\\
\mathbf{E}\left[  \sup_{T_{n\alpha}}X(k,\tilde{t},\tilde{x})\right]   &  \leq
Ktd\sqrt{n\alpha}, \label{stimasup}%
\end{align}
where $K$ is a constant that depends on the covariance of the environment $Q$.
Inequality (\ref{stimaprob}) can be found textually in \cite{FV}. Inequality
(\ref{stimasup}) is established identically to equation (30) in \cite{FV},
with the minor difference that the total number of paths in $\mathcal{P}_{m}$
is not $2^{m}$ but $\left(  2d\right)  ^{m}$, which, in the inequality above
(30) near the bottom of page 33 in \cite{FV}, accounts for a factor
$e^{1+\log(6d)}=6ed$ instead of $e^{c_{1}}$ therein, hence the factor $d$ in
(\ref{stimasup}).

Defining $Y_{n\alpha}=\sup_{T_{n\alpha}}X(k,\tilde{t},\tilde{x})$, we can now
bound $p_{t}^{\varepsilon}(\beta)$ as follows:
\[
tp_{t}^{\varepsilon}(\beta)\leq\mathbf{E}\left[  \log(A+B)\right]
\leq\mathbf{E}\left[  \left(  \log A\right)  _{+}\right]  +\mathbf{E}\left[
\left(  \log B\right)  _{+}\right]  +\log2,
\]
where%
\[
A=P_{b}\left[  |N_{t}|\leq\alpha t\right]  \exp\left(  2\beta Y_{\alpha
}\right)  ,\quad\mbox{and}\quad B=\sum_{n\geq1}P_{b}\left[  |N_{t}|\in\left[
n\alpha t,(n+1)\alpha t\right]  \right]  \exp\left(  2\beta Y_{\alpha
(n+1)}\right)  .
\]
We will now bound the terms $A$ and $B$ separately.\vspace{0.3cm}

\noindent\textit{Step 4: The factor A.} We can bound $P_{b}\left[  |N_{t}%
|\leq\alpha t\right]  $ by 1 and we easily get, invoking (\ref{stimasup}),
\begin{equation}
\label{A}\mathbf{E}\left[  \left(  \log A\right)  _{+}\right]  \leq
2\beta\mathbf{E}\left[  Y_{\alpha}\right]  \leq2\beta K d t\sqrt\alpha.
\end{equation}

\noindent\textit{Step 5: The factor B.} Let $\mu=\mathbf{E}\left[
Y_{\alpha(n+1)}\right]  $. Since $X$ is a Gaussian field and since it is easy
to show that
\[
\sigma^{2}:=\sup_{(m,\tilde{t},\tilde{x})}\mbox{\textbf{Var}}(X(k,\tilde
{t},\tilde{x}))\leq tQ(0),
\]
the so called Borell-Sudakov inequality (see \cite{A} or \cite{VV}) implies
that, for a constant $a>0$,%

\begin{equation}
\mathbf{E}\left[  \exp\left(  a\left\vert Y_{\alpha n}-\mu\right\vert \right)
\right]  \leq2\exp\left(  \frac{a^{2}\sigma^{2}}{2}\right)  =2\exp\left(
\frac{a^{2}tQ(0)}{2}\right)  .\label{BS}%
\end{equation}
Fix now a number $\gamma\in\left(  1/2,1\right)  $ and let us denote $\log
_{+}(B)=(\log B)_{+}$ . We have
\begin{multline*}
\frac{1}{t^{\gamma}}\mathbf{E}\left[  \log_{+}B\right]  =\mathbf{E}\left[
\log_{+}\left(  \sum_{n\geq1}P_{b}\left[  |N_{t}|\in\left[  nt\alpha
,(n+1)t\alpha\right]  \right]  \exp\left(  2\beta Y_{\alpha(n+1)}\right)
\right)  ^{t^{-\gamma}}\right]  \\
\leq\mathbf{E}\left[  \log_{+}\left(  \sum_{n\geq1}P_{b}\left[  |N_{t}%
|>nt\alpha\right]  \exp\left(  2\beta(Y_{\alpha(n+1)}-\mu)\right)  \exp\left(
2\beta Ktd\sqrt{\alpha(n+1)}\right)  \right)  ^{t^{-\gamma}}\right]  ,
\end{multline*}
where we used that (\ref{stimasup}) implies $\mu\leq Kdt\sqrt{(n+1)\alpha}$.
We also know that for any sequence of non-negative reals $(x_{n})_{n}$ the
following holds: $(\sum_{n}x_{n})^{t^{-\gamma}}\leq\sum_{n}x_{n}^{t^{-\gamma}%
}$. Thus we have
\begin{align*}
&  \frac{1}{t^{\gamma}}\mathbf{E}\left[  \log_{+}B\right]  \\
&  \leq\mathbf{E}\left[  \log_{+}\left(  \sum_{n\geq1}\left(  P_{b}\left[
|N_{t}|>nt\alpha\right]  \right)  ^{t^{-\gamma}}\exp\left(  \frac{2\beta
}{t^{\gamma}}(Y_{\alpha(n+1)}-\mu)\right)  \exp\left(  2t^{1-\gamma}\beta
Kd\sqrt{\alpha(n+1)}\right)  \right)  \right]  \\
&  \leq\mathbf{E}\left[  \log_{+}\left[  d^{t^{-\gamma}}\sum_{n\geq1}%
\exp\left(  \frac{2\beta}{t^{\gamma}}(Y_{\alpha(n+1)}-\mu)\right)  \exp\left(
-\frac{t^{1-\gamma}}{2}y_{n}\right)  \right]  \right]  ,
\end{align*}
where we used estimate (\ref{stimaprob}) in the following way:
\begin{align*}
P_{b}\left[  |N_{t}|>nt\alpha\right]   &  \leq\sum_{j=1}^{d}P_{b}\left[
N_{t}^{j}>\frac{nt\alpha}{d}\right]  =dP_{b}\left[  N_{t}^{1}>\frac{nt\alpha
}{d}\right]  \\
&  \leq d\exp\left(  -\frac{t}{2}\left(  \frac{\alpha n\varepsilon}{d}\right)
^{2}+\frac{t\alpha n}{d}\right)  ,
\end{align*}
and where we have obtained:
\[
y_{n}=\left(  \frac{\varepsilon\alpha n}{d}\right)  ^{2}-\frac{2\alpha n}%
{d}-4\beta Kd\sqrt{\alpha(n+1)}.
\]
Now, bounding $\log_{+}(x)$ from above by $\log(1+x)$, for $x\geq1$, and using
Jensen's inequality, we have:
\[
\frac{1}{t^{\gamma}}\mathbf{E}\left[  \log_{+}B\right]  \leq\log\left[
1+\sum_{n\geq1}\mathbf{E}\left[  \exp\left(  \frac{2\beta}{t^{\gamma}%
}(Y_{\alpha(n+1)}-\mu)\right)  \right]  \exp\left(  \frac{-t^{1-\gamma}}%
{2}y_{n}\right)  \right]  ,
\]
so, using (\ref{BS}), it is readily checked that
\[
\frac{1}{t^{\gamma}}\mathbf{E}\left[  \log_{+}B\right]  \leq\log\left[
1+2\exp\left(  \frac{2\beta^{2}Q(0)}{t^{2\gamma-1}}\right)  \sum_{n\geq1}%
\exp\left(  \frac{-t^{1-\gamma}}{2}y_{n}\right)  \right]  .
\]
In order for the series above to converge, we must choose $\alpha$ so as to
compensate the negative terms in $y_{n}$. Specifically, we choose
\begin{equation}
\left(  \frac{\alpha\varepsilon}{d}\right)  ^{2}=16\beta Kd\sqrt{\alpha}%
,\quad\mbox{ i.e. }\quad\alpha=(16\beta Kd^{3}\varepsilon^{-2})^{2/3}%
.\label{alpha}%
\end{equation}
With this choice, we end up with:
\[
y_{n}=\left(  \frac{\alpha\varepsilon}{d}\right)  ^{2}\left(  n^{2}-\frac
{2dn}{\alpha\varepsilon^{2}}-\frac{1}{4}\sqrt{n+1}\right)  .
\]
Now we note that:
\begin{equation}
\mbox{If we choose $\varepsilon,\beta$ such that $\beta\varepsilon\geq
d^{-3/2}$}\quad\Rightarrow\quad\frac{\alpha\varepsilon^{2}}{d}=\left(
16K\beta\varepsilon\right)  ^{2/3}d\geq4,\label{eq:cdt-beta-ep}%
\end{equation}
so that
\[
y_{n}\geq\left(  \frac{\alpha\varepsilon}{d}\right)  ^{2}\left(  n^{2}%
-\frac{n}{2}-\frac{1}{4}\sqrt{n+1}\right)  ,
\]
and since $n^{2}-\frac{n}{2}-\frac{\sqrt{n+1}}{4}\geq\frac{n}{8}$, we get
\begin{align*}
&  \sum_{n\geq1}\exp\left(  -\frac{t^{1-\gamma}}{2}\left(  \frac
{\alpha\varepsilon}{d}\right)  ^{2}\left(  n^{2}-\frac{2dn}{\alpha
\varepsilon^{2}}-\frac{1}{4}\sqrt{n+1}\right)  \right)  \\
&  \leq\sum_{n\geq1}\exp\left(  -\frac{t^{1-\gamma}}{2}\left(  \frac
{\alpha\varepsilon}{d}\right)  ^{2}\frac{n}{8}\right)  =\frac{1}{1-\exp\left(
-\frac{1}{16}t^{1-\gamma}\left(  \frac{\alpha\varepsilon}{d}\right)
^{2}\right)  }-1.
\end{align*}
Notice that this last term can be made smaller than $1$ if $t$ is large
enough. Hence we can write a final estimate on $\mathbf{E}\left[  \log
_{+}B\right]  $ as follows: for large $t$ we have
\begin{align}
\frac{1}{t^{\gamma}}\mathbf{E}\left[  \log_{+}B\right]   &  \leq\log\left[
1+2d^{t^{-\gamma}}\exp\left(  \frac{2\beta^{2}Q(0)}{t^{2\gamma-1}}\right)
\right]  \nonumber\\
&  \leq\log(1+2d^{t^{-\gamma}})+\frac{2\beta^{2}Q(0)}{t^{2\gamma-1}}.\label{B}%
\end{align}

\noindent\textit{Final step.} Using inequalities (\ref{A}) and (\ref{B}) and
the value of $\alpha$, we can estimate $p_{t}^{\varepsilon}(\beta)$ in the
following way:
\begin{align*}
p_{t}^{\varepsilon}(\beta)  &  \leq2\beta Kd\sqrt{\alpha}+\frac{\log2}%
{t}+\frac{\log(1+2d^{t^{-\gamma}})}{t^{1-\gamma}}+\frac{2\beta^{2}%
Q(0)}{t^{\gamma}}\\
&  \leq2\beta Kd\sqrt{\alpha}+o(1).
\end{align*}
So using the value of $\alpha$ given in (\ref{alpha}) we have
\begin{equation}
p_{t}^{\varepsilon}(\beta)\leq C\frac{\beta^{4/3}d^{2}}{\varepsilon^{2/3}%
}+o(1), \label{eq:low-bnd-p-t-ep}%
\end{equation}
where $C$ is a constant that depends on $Q$ and that can change from line to
line. Putting this result in (\ref{pt}) and taking the limit for
$t\rightarrow\infty$ we get
\[
\limsup_{t\rightarrow\infty}p_{t}(\beta)\leq C\left(  \beta^{2}d^{H}%
\varepsilon^{2H}+d^{2}\beta^{4/3}\varepsilon^{-2/3}\right)  .
\]
In order to make this upper bound as small as possible we can choose
$\varepsilon$ such that
\[
\beta^{2}d^{H}\varepsilon^{2H}=d^{2}\beta^{4/3}\varepsilon^{-2/3}%
,\quad\mbox{ i.e. }\quad\varepsilon=d^{\frac{6-3H}{2+6H}}\beta^{-\frac
{1}{1+3H}},
\]
so that
\[
\limsup_{t\rightarrow\infty}p_{t}(\beta)\leq C\beta^{\frac{2+4H}{1+3H}%
}d^{\frac{7H}{1+3H}},
\]
which is the announced result. We then only need to check for what values of
$\beta$ we are allowed to make this choice of $\varepsilon$. Condition
(\ref{hub}) states that we must use $\varepsilon\leq r_{1}$. This is
equivalent to $\beta\geq\beta_{0}=:\left(  r_{1}\right)  ^{-1-3H}d^{3-3H/2}$.
One can check in this case that the restriction on $\varepsilon,\beta$ in
(\ref{eq:cdt-beta-ep}) is trivially satisfied.
\end{proof}

\subsection{Lower bound in the Brownian case}

\label{sec:low-bnd-brownian}

In the following proposition, which implies the lower bound in Theorem
\ref{thm:bnd-free-nrj-brown-holder}, we shall also try to specify the
dependence of the constants with respect to the dimension $d$. Let us state an
interesting feature of this dependence. The proof of the proposition below
shows that the results it states hold only for $\beta\geq\beta_{0}=cd^{1-H/2}%
$. One may ask the question of what happens to the behavior of the partition
function when the dimension is linked to the inverse temperature via the
relation $\beta=\beta_{0}$, and one allows the dimension to be very large. The
lower bounds on the value $p\left(  \beta\right)  $ in the proposition below
will then increase, and while they must still not exceed the global bound
$\beta^{2}Q\left(  0\right)  /2$, the behavior for large $\beta$ turns out to
be quadratic in many cases. The reader will check that, when $H>1/2$, this
translates as $p\left(  \beta\right)  \geq c\beta^{2/(2-H)}$ which is
quadratic when $H=1$, and $p\left(  \beta\right)  \geq c\beta^{2}$ for all
$H\leq1/2$. This is an indication that for extremely high dimensions and
inverse temperatures, for $H\leq1/2$ or $H=1$, strong disorder may not hold.
Strong disorder for Brownian polymers may break down for complex,
infinite-dimensional polymers. This is only tangential to our presentation, however.

\begin{proposition}
\label{prop:lower-bnd-holder} Recall that $\delta$ has been defined at
(\ref{eq:def-delta}) and assume that there exist a number $H\in\left(
0,1\right]  $ and some positive constants $c_{2}$, $r_{2}$ such that for all
$x,y\in\mathbb{R}^{d}$ with $\left\vert x-y\right\vert \leq r_{2}$, we have
\begin{equation}
\delta(x-y)>c_{2}\left\vert x-y\right\vert ^{H}. \label{hlb}%
\end{equation}
Then if $H\leq1/2$, there exists a constant $C$ depending only on $Q$, and a
constant $\beta_{0}$ depending only on $Q$ and $d$, such that, for all
$\beta>\beta_{0}$,%
\[
p(\beta)\geq Cd^{\frac{2H-1}{H+1}}\beta^{\frac{2}{H+1}}.
\]
On the other hand if $H>1/2$, there exists a constant $C^{\prime}$ depending
only on $Q$, and a constant $\beta_{0}^{\prime}$ depending only on $Q$ and
$d$, such that for all $\beta>\beta_{0}^{\prime}$
\[
p(\beta)\geq C^{\prime}d^{\frac{2H-1}{3}}\beta^{\frac{4}{3}}.
\]

\end{proposition}

\begin{proof}
Here again, we divide the proof in several steps.\vspace{0.3cm}

\noindent\textit{Step 1: Strategy.} From relation (\ref{eq:lim-p-t-beta}), we
have
\[
p(\beta)=\sup_{t\geq0}p_{t}(\beta),
\]
where $p_{t}(\beta)$ is defined by equation (\ref{fren}). So a lower bound for
$p(\beta)$ will be obtained by evaluating $p_{t}(\beta)$ for any fixed value
$t$. Additionally, by the positivity of the exponential factor in the
definition of $Z_{t}$, one may include as a factor inside the expectation
$E_{b}$ the sum of the indicator functions of any disjoint family of events of
$\Omega_{b}$. In fact, we will need only two events, which will give the main
contribution to $Z_{t}$ at a logarithmic scale.\vspace{0.3cm}

\noindent\textit{Step 2: Setup.} Let $A_{+}(b)$ and $A_{-}(b)$ be two
disjoints events defined on the probability space $\Omega_{b}$ under $P_{b}$,
which will be specified later on. Set
\[
X_{b}=-\beta H_{2t}=\beta\int_{0}^{2t}W(ds,b_{s}).
\]
Conditioning by the two events $A_{+}(b)$ and $A_{-}(b)$ and using Jensen's
inequality we have
\begin{equation}
\mathbf{E}(\log Z_{t})\geq\log\left(  \min\left\{  P_{b}(A_{+});P_{b}%
(A_{-})\right\}  \right)  +\mathbf{E}\left[  \max\left\{  \tilde{Z}_{+}%
;\tilde{Z}_{-}\right\}  \right]  ,\label{p2t}%
\end{equation}
where
\[
\tilde{Z}_{+}:=E_{b}\left[  X_{b}\mid A_{+}\right]  \quad\mbox{ and }\quad
\tilde{Z}_{-}:=E_{b}\left[  X_{b}\mid A_{-}\right]  .
\]
These two random variables form a pair of centered jointly Gaussian random
variables: indeed they are both limits of linear combinations of values of a
single centered Gaussian field. Thus this implies
\[
\mathbf{E}\left[  \max\left\{  \tilde{Z}_{+};\tilde{Z}_{-}\right\}  \right]
=\frac{1}{\sqrt{2\pi}}\left(  \mathbf{E}\left[  \left(  \tilde{Z}_{+}%
-\tilde{Z}_{-}\right)  ^{2}\right]  \right)  ^{1/2}.
\]
Therefore we only have to choose sets $A_{+}$ and $A_{-}$ not too small, but
still decorrelated enough so that condition (\ref{hlb}) guarantees a certain
amount of positivity in the variance of $\tilde{Z}_{+}-\tilde{Z}_{-}$%
.\vspace{0.3cm}

\noindent\textit{Step 3: Choice of }$A_{+}$\textit{ and }$A_{-}$\textit{.} Let
$f$ be a positive increasing function. We take
\[
A_{+}=\left\{  f(t)\leq b_{s}^{i}\leq2f(t),\forall i=1,\ldots,d,\quad\forall
s\in\left[  t,2t\right]  \right\}  ,
\]
\[
A_{-}=\left\{  -2f(t)\leq b_{s}^{i}\leq-f(t),\forall i=1,\ldots,d,\quad\forall
s\in\left[  t,2t\right]  \right\}  .
\]
In other words, we force each component of our trajectory $b$ to be, during
the entire time interval $[t,2t]$, in one of two boxes of edge size $f\left(
t\right)  $ which are at a distance of $2f\left(  t\right)  $ from each other.
Because these two boxes are symmetric about the starting point of $b,$ the
corresponding events have the same probability. While this probability can be
calculated in an arguably explicit way, we give here a simple lower bound
argument for it. Using time scaling, the Markov property of Brownian motion,
the notation $a=f\left(  t\right)  /\sqrt{t}$, we have
\begin{align}
P_{b}(A_{+}) &  =\prod_{i=1}^{d}P_{b}\left(  \forall s\in\left[  1,2\right]
:b_{s}^{i}\in\left[  a,2a\right]  \right)  \nonumber\\
&  =\prod_{i=1}^{d}\frac{1}{2\pi}\int_{a}^{2a}P_{b}\left(  \forall s\in\left[
0,1\right]  :b_{s}^{i}+y\in\left[  a,2a\right]  \right)  e^{-y^{2}%
/2}dy\nonumber\\
&  \geq\prod_{i=1}^{d}\frac{1}{2\pi}\int_{5a/4}^{7a/4}P_{b}\left(  \forall
s\in\left[  0,1\right]  :b_{s}^{i}+y\in\left[  y-\frac{a}{4},y+\frac{a}%
{4}\right]  \right)  e^{-y^{2}/2}dy\nonumber\\
&  =\left[  P_{b}\left(  b_{1}^{1}\in\left[  5a/4,7a/4\right]  \right)
P_{b}\left(  \forall s\in\left[  0,1\right]  :|b_{s}^{1}|\leq a/4\right)
\right]  ^{d}.\label{PA}%
\end{align}

\noindent\textit{Step 4: Estimation of} $\tilde{Z}_{+}$ \textit{and}
$\tilde{Z}_{-}$\textit{.} It was established in \cite{FV} that in dimension
$d=1$
\[
\mathbf{E}\left[  \left(  \tilde{Z}_{+}-\tilde{Z}_{-}\right)  ^{2}\right]
\geq\beta^{2}\int_{t}^{2t}\mathbf{E}\left[  \left(  \delta(x_{s,+}^{\ast
}-x_{s,-}^{\ast})\right)  ^{2}\right]  ds
\]
where the quantities $x_{s,+}^{\ast}$ and $x_{s,-}^{\ast}$ are random
variables such that for all $s\in\left[  t,2t\right]  $: $x_{s,+}^{\ast}%
\in\left[  f(t),2f(t)\right]  $ and $x_{s,-}^{\ast}\in\left[
-2f(t),-f(t)\right]  $. In dimension $d\geq1$ the result still holds. In fact
in this case we have $x_{s,+}^{\ast},x_{s,-}^{\ast}\in\mathbb{R}^{d}$, so it
is sufficient to take each component of the $x_{s,+}^{\ast}$ in the interval
$\left[  f(t),2f(t)\right]  $ and each component of $x_{s,-}^{\ast}$ in
$\left[  -2f(t),-f(t)\right]  $, so their distance is greater than
$d^{1/2}f(t)$. Thus, using condition (\ref{hlb}), we have
\begin{equation}
\mathbf{E}\left[  \left(  \tilde{Z}_{+}-\tilde{Z}_{-}\right)  ^{2}\right]
\geq\beta^{2}\int_{t}^{2t}C\left\vert x_{s,+}^{\ast}-x_{s,-}^{\ast}\right\vert
^{2H}ds\geq Ct\beta^{2}d^{H}\left(  f(t)\right)  ^{2H}%
,\label{eq:expec-sqr-z+-z}%
\end{equation}
where as usual $C$ is a constant that can change from line to line. Hence, we
obtain:
\begin{equation}
\mathbf{E}\left[  \max\left\{  \tilde{Z}_{+};\tilde{Z}_{-}\right\}  \right]
=\frac{1}{\sqrt{2\pi}}\left(  \mathbf{E}\left[  \left(  \tilde{Z}_{+}%
-\tilde{Z}_{-}\right)  ^{2}\right]  \right)  ^{1/2}\geq C\beta\sqrt{t}\left(
f(t)\right)  ^{H}d^{H/2},\label{E}%
\end{equation}
Observe that in order to use condition (\ref{hlb}) we have to impose $f(t)\leq
r_{2}$.\vspace{0.3cm}

\noindent\textit{Step 5: The case }$H\leq1/2$\textit{.} It is possible to
prove that in this case the optimal choice for $f$ is $f(t)=\sqrt{t}$, which
corresponds to $a=1$, so that $P_{b}(A_{+})$ is a universal constant that does
not depend on $t$. Thus we have, from (\ref{p2t}), (\ref{PA}) and (\ref{E}),
for any $t>0$,
\begin{equation}
p_{2t}(\beta)=\frac{\mathbf{E}\left[  \log Z_{2t}\right]  }{2t}\geq\frac{d\log
C}{2t}+C\beta d^{H/2}t^{\frac{H-1}{2}}.\label{eq:low-bnd-p-2t}%
\end{equation}
Now we may maximize the above function over all possible values of $t>0$. To
make things simple, we choose $t$ so that the second term equals twice the
first, yielding $t$ of the form $t=Cd^{\frac{2-H}{H+1}}\beta^{-\frac{2}{H+1}%
},$ and therefore
\[
\sup_{t>0}p_{2t}(\beta)\geq Cd^{\frac{2H-1}{H+1}}\beta^{\frac{2}{H+1}}.
\]
This result holds as long as the use of condition (\ref{hlb}) can be
justified, namely as long as $f(t)\leq r_{2}$. This is achieved as soon as
$\beta>\beta_{0}$ where $\beta_{0}=Cr_{2}^{-H-1}d^{1-H/2},$ and since
$H\leq1/2$, $\beta_{0}\geq Cd^{3/4}$.\vspace{0.3cm}

\noindent\textit{Step 6: The case }$H>1/2$\textit{.} In this case we consider
$f(t)=ct^{\alpha}$, for a given $\alpha\in\lbrack0,1/2)$ and some constant $c$
chosen below. Thus we have $a=ct^{\alpha-1/2}$. In this case, if $a$ is larger
than a universal constant $K_{u}$, the result (\ref{PA}) yields that, for some
constant $C$, we have
\[
P_{b}(A_{+})\geq\prod_{i=1}^{d}\exp(-Ca^{2})=\exp(-Cc^{2}dt^{2\alpha-1}).
\]
So, using again condition (\ref{hlb}) and relation (\ref{E}) we obtain
\[
p_{2t}(\beta)\geq-Cdt^{2\alpha-2}+C\beta d^{H/2}t^{\alpha H-1/2},
\]
where the constant $C$ may also include the factor $c^{2}$. Again, choosing
$t$ so that the second term equals twice the first, we have
\begin{equation}
t=Cd^{\frac{1-H/2}{\alpha(H-2)+3/2}}\beta^{-\frac{1}{\alpha(H-2)+3/2}},
\label{eq:value-t-lower-bnd-holder}%
\end{equation}
and so
\[
\sup_{t>0}p_{2t}(\beta)\geq Cd^{\frac{H-1/2}{\alpha(H-2)+3/2}}\beta
^{-\frac{2\alpha-2}{\alpha(H-2)+3/2}}.
\]
In order to maximize the power of $\beta$ in the lower bound for $\sup
_{t>0}p_{t}(\beta)$ we should find the maximum of the function
\[
g(\alpha)=\frac{2-2\alpha}{\alpha(H-2)+3/2}%
\]
for $0\leq\alpha<1/2.$ Since this function is monotone decreasing when
$H>1/2$, the maximum is reached for $\alpha=0$, so $g(0)=4/3$.

Recall once again that, in order to apply condition (\ref{hlb}) in the
computations above, we had to assume $f(t)\leq r_{2}$; since now $f\left(
t\right)  $ is the constant $c$, we only need to choose $c=r_{2}$. We also had
to impose $a=r_{2}t^{-1/2}>K_{u}$, which translates as $\beta>\beta
_{0}^{\prime}:=\left(  K_{u}/r_{2}\right)  ^{4/3}d^{1-H/2}$.
\end{proof}

\subsection{Logarithmic regularity scale}

As mentioned in the introduction, the special shape of our Figure
\ref{fig:exp-beta} induces us to explore the regions of low spatial regularity
for $W$, in order to investigate some new possible scaling in the strong
disorder regime. In other words, we shall work in this section under the
assumptions that there exist positive constants $c_{0}$, $c_{1}$, and $r_{1}$,
and $\beta\in(0,\infty)$, such that for all $x,y$ with $\left\vert
x-y\right\vert \leq r_{1}$,%
\begin{equation}
c_{0}\log^{-\gamma}\left(  1/\left\vert x-y\right\vert \right)  \leq
\delta\left(  x-y\right)  \leq c_{1}\log^{-\gamma}\left(  1/\left\vert
x-y\right\vert \right)  , \label{loghypo2}%
\end{equation}
where $\gamma>0$. Assumption (\ref{loghypo2}) implies that $W$ is not
spatially H\"{o}lder-continuous for any exponent $H\in(0,1]$. Moreover, the
theory of Gaussian regularity implies that, if $\gamma>1/2$, $W$ is
almost-surely continuous in space, with modulus of continuity proportional to
$\log^{-\gamma+1/2}\left(  1/\left\vert x-y\right\vert \right)  $, while if
$\gamma\leq1/2$, $W$ is almost-surely not uniformly continuous on any interval
in space. The case $\gamma=1/2$, which is the threshold between continuous and
discontinuous $W$, is of special interest, since it can be related to the
discrete space polymer which will be studied in the next section. The main
result which will be proved here is the following:

\begin{theorem}
\label{thmlog}Assume condition (\ref{loghypo2}). We have for some constants
$C_{0}$ and $C_{1}$ depending only on $Q$, for all $\beta$ large enough,
\[
C_{0}\frac{\beta^{2}}{d}\log^{-2\gamma}\left(  \frac{\beta}{\sqrt d}\right)
\leq p(\beta) \leq C_{1}\beta^{2}\log^{-2\gamma}\left(  \frac{\beta}{\sqrt
d}\right)  .
\]

\end{theorem}

\begin{proof}
\noindent\emph{Step 1: }\emph{Setup.} Nearly all the calculations in the proof
of Propositions \ref{prop:upper-bnd-holder} and \ref{prop:lower-bnd-holder}
are still valid in our situation.\vspace*{0.1in}

\noindent\emph{Step 2: }\emph{Lower bound.} For the lower bound, reworking the
argument in Step 2 in the proof of Proposition \ref{prop:lower-bnd-holder},
using the function $\log^{-\gamma}\left(  x^{-1}\right)  $ instead of the
function $x^{H}$, we obtain the following instead of (\ref{eq:expec-sqr-z+-z}%
):
\[
\mathbf{E}\left[  \left(  Z_{+}-Z_{-}\right)  ^{2}\right]  \geq t\left(  \beta
c_{0}\right)  ^{2}\left(  \log\left(  \frac{1}{\sqrt{d}f(t)}\right)  \right)
^{-2\gamma},
\]
which implies, instead of (\ref{eq:low-bnd-p-2t}) in Step 5 of that proof, the
following:%
\[
p_{2t}(\beta)\geq\frac{d\log C}{2t}+C\beta t^{-1/2}\left(  \log\left(
\frac{1}{\sqrt{d}f(t)}\right)  \right)  ^{-\gamma}.
\]
In other words, now choosing $f\left(  t\right)  =t^{1/2}$ as we did in the
case $H<1/2$ (recall that we are in the case of small $H$, as stated in the
introduction),
\[
p_{2t}(\beta)\geq\frac{d\log C}{2t}+C\beta t^{-1/2}\left(  \log\left(
\frac{1}{\sqrt{dt}}\right)  \right)  ^{-\gamma}.
\]
Now choose $t$ such that the second term in the right-hand side above equals
twice the first, i.e.
\[
t^{1/2}\log^{-\gamma}\left(  \frac{1}{\sqrt{dt}}\right)  =Cd\beta^{-1}.
\]
For small $t$, the function on the left-hand side is increasing, so that the
above $t$ is uniquely defined when $\beta$ is large. We see in particular that
when $\beta$ is large, $t$ is small, and we have $t^{-1}\leq\beta^{2}$. This
fact is then used to imply%
\[
\frac{1}{t}=\left(  \frac{C\beta}{d}\right)  ^{2}\left(  \log\left(  \frac
{1}{\sqrt{dt}}\right)  \right)  ^{-2\gamma}\geq2\left(  C\beta\right)
^{2}\log^{-2\gamma}\left(  \beta\right)  .
\]
Therefore, for some constants $\beta_{2}$ and $c$ depending only on $Q$, for
the $t$ chosen above with $\beta\geq\beta_{2}$,%
\[
p_{2t}(\beta)\geq\frac{C\beta^{2}}{d}\left(  \log\left(  \frac{\beta}{\sqrt
{d}}\right)  \right)  ^{-2\gamma}.
\]
\vspace*{0.1in}

\noindent\emph{Step 3: }\emph{Upper bound.} Here, returning to the proof of
Proposition \ref{prop:upper-bnd-holder}, the upper bound
(\ref{eq:low-bnd-p-t-ep}) in the final step of that proof holds regardless of
$\delta$, and therefore, using the result of Remark
\ref{rmk:bdn-error-pt-pt-ep} with $\delta\left(  r\right)  =\log^{-\gamma
}\left(  1/r\right)  $, we immediately get that there exists $c$ depending
only on $Q$ such that for all $\varepsilon<r_{1}$ and all $\beta>\beta_{3}$,%
\[
\limsup_{t\rightarrow\infty}p_{t}(\beta) \leq C\beta^{2} \left(  \log\left(
\frac{1}{\varepsilon\sqrt d}\right)  \right)  ^{-2\gamma} +C d^{2}\beta
^{4/3}\varepsilon^{-2/3},
\]
as long as one is able to choose $\varepsilon$ so that $\beta\varepsilon\geq
1$. By equating the two terms in the right-hand side of the last inequality
above, we get%
\[
\varepsilon\left(  \log\left(  \frac{1}{\varepsilon\sqrt d}\right)  \right)
^{-3\gamma} =C d^{3}\beta^{-1}.
\]
Since the function $\varepsilon\mapsto\varepsilon\log^{-3\gamma}\left(
1/(\varepsilon\sqrt{d})\right)  $ is increasing for small $\varepsilon$, the
above equation defines $\varepsilon$ uniquely when $\beta$ is large, and in
that case $\varepsilon$ is small. We also see that for any $\theta>0$, for
large $\beta$, $1/\varepsilon\geq\beta^{1-\theta}$. Therefore we can write,
for $\beta\geq\beta_{3}$, almost surely,%
\[
\limsup_{t\rightarrow\infty}p_{t}(\beta)\leq C \left(  1-\theta\right)
^{-2\gamma} \beta^{2} \left(  \log\left(  \frac{\beta}{\sqrt d}\right)
\right)  ^{-2\gamma}.
\]
This finishes the proof of the theorem.
\end{proof}

\section{Estimates of the free energy: discrete space}

Recall that, up to now, we have obtained our bounds on the free energy in the
following manner: the upper bound has been computed by evaluation of the
supremum of a well-chosen random Gaussian field, while the lower bound has
been obtained by introducing two different events, depending on the Brownian
configuration, which capture most of the logarithmic weight of our polymer
distribution. This strategy also works in the case of the random walk polymer
whose Hamiltonian is described by (\ref{eq:expression-hamil-discrete}),
without many additional efforts, but a separate proof is still necessary. This
section shows how this procedure works, resulting in the proof of Theorem
\ref{thm:bnd-free-nrj-walk}.

\vspace{0.3cm}

Quantities referring to the random walk polymer have been denoted by $\hat
{b},\hat{W},\hat{E}_{\hat{b}},\hat{\mbox{{\bf E}}}$, etc... In this section,
for notational sake, we will omit the hats in the expressions above, and write
instead $b,W,E_{b},\mbox{{\bf E}}$ like in the Brownian case. Recall our
simple non-degeneracy condition on $Q$ in this case:
\begin{equation}
c_{Q}:=\sup_{1\leq i\leq d}\left(  Q(0)-Q(2e_{i})\right)  ^{1/2}>0,
\label{eq:cdt-q-discrete}%
\end{equation}
where $e_{i},\,i=1,\cdots,d$ are the unit vectors in ${\mathbb{Z}}^{d}$.
Condition (\ref{eq:cdt-q-discrete}), which is used only in the lower bound
result, is extremely weak. It essentially covers all possible homogeneous
covariance functions, except the trivial one $Q\left(  x\right)  \equiv
Q\left(  0\right)  $ for all $x$, which is the case where $W$ does not depend
on $x$, in which case the Hamiltonian has no effect. Indeed, assume that there
exists an $x_{0}\in{\mathbb{Z}}^{d}$ such that $W\left(  t,0\right)  $ and
$W\left(  t,x_{0}\right)  $ are not (a.s.) equal. Then $Q\left(  x_{0}\right)
<Q\left(  0\right)  $. Our lower bound proof below can then be adapted to use
this condition instead of Condition (\ref{eq:cdt-q-discrete}). We do not
comment on this point further.

\subsection{Lower bound for the random walk polymer}

The lower bound announced in Theorem \ref{thm:bnd-free-nrj-walk} is contained
in the following.

\begin{proposition}
Assume condition (\ref{eq:cdt-q-discrete}) holds true. Then there exists a
constant $\beta_{0}>0$, which depends on $d$ and on $c_{Q}$ and a constant
$C>0$, which depend only on $c_{Q}$, such that if $\beta>\beta_{0}$ then
almost surely
\[
\lim_{t\rightarrow\infty}\frac{1}{t}\log Z_{t}\geq C\frac{\beta^{2}}{\log
\beta}.
\]

\end{proposition}

\begin{proof}
Invoking Proposition \ref{prop:simple-prop-free-discrete}, we have $p\left(
\beta\right)  =\lim_{t\rightarrow\infty}p_{t}(\beta)=\sup_{t\geq0}p_{t}%
(\beta)$. Therefore, any chosen fixed value $t$ yields $p_{t}(\beta)$ as a
lower bound for $p\left(  \beta\right)  $.

We express $Z_{t}$ by using the fact that each component of $b$ is constant
between its jump times, which are uniformly distributed on the simplex, given
$N_{t}$ the total number of jumps before time $t$, which is a Poisson r.v.
with parameter $2dt$. Moreover the visited sites $\left(  x_{k}\right)
_{k=1}^{N_{t}}$ are uniformly distributed on the set of all nearest-neighbor
paths of length $N_{t}$ started at $0$, given $N_{t}$. For a lower bound on
$p_{t}(\beta)$, we throw out, in the expectation defining $Z_{t}$, all the
paths $b$ that do not jump exactly once before time $t$. We also throw out all
jump positions that are not $\pm e_{i}$, where $c_{Q}=(Q(0)-Q(2e_{i}%
))^{1/2}>0$. Therefore,
\[
Z_{t}\geq P_{b}\left[  N_{t}=1\right]  \frac{1}{2d}\int_{0}^{t}\frac{ds}%
{t}\left(  e^{\beta W\left(  s,0\right)  +\beta W\left(  [s,t],e_{i}\right)
}+e^{\beta W\left(  s,0\right)  +\beta W\left(  [s,t],-e_{i}\right)  }\right)
,
\]
where $W\left(  [s,t],x\right)  :=W\left(  t,x\right)  -W\left(  s,x\right)
$. Here, given $N_{t}=1$, $1/(2d)$ is the weight of the path that jumps to
$\pm e_{i}$, and $\mathbf{1}_{[0,t]}\left(  s\right)  ds/t$ is the law of the
single jump time. Using this and Jensen's inequality, we get%
\[
Z_{t}\geq dte^{-2td}\int_{0}^{t}\frac{ds}{t}\left(  e^{\beta W\left(
s,0\right)  +\beta W\left(  [s,t],e_{i}\right)  }+e^{\beta W\left(
s,0\right)  +\beta W\left(  [s,t],-e_{i}\right)  }\right)  ,
\]%
\[
\frac{1}{t}\mathbf{E}(\log Z_{t})\geq\frac{\log t}{t}-2d+\beta\int_{0}%
^{t}\frac{ds}{t^{2}}\mathbf{E}\left[  \max\left(  W([s,t],e_{i}%
);W([s,t],-e_{i})\right)  \right]  .
\]

Now we evaluate the expected maximum above. The vector $\left(  W\left(
[s,t],e_{i}\right)  ,W\left(  [s,t],-e_{i}\right)  \right)  $ is jointly
Gaussian with common variances $\sqrt{t-s}Q(0)$ and covariance $\sqrt
{t-s}Q(2)$. Therefore%
\[
\mathbf{E}\left[  \max\left(  W\left(  [s,t],e_{i}\right)  ,W\left(
[s,t],-e_{i}\right)  \right)  \right]  =\frac{1}{2}\mathbf{E}\left[
\left\vert W\left(  [s,t],e_{i}\right)  -W\left(  [s,t],-e_{i}\right)
\right\vert \right]
\]%
\begin{equation}
=\frac{1}{\sqrt{2\pi}}\left(  \mbox{\bf Var}\left[  W\left(  [s,t],e_{i}%
\right)  -W\left(  [s,t],-e_{i}\right)  \right]  \right)  ^{1/2}=\frac
{1}{\sqrt{\pi}}\sqrt{t-s}\sqrt{Q\left(  0\right)  -Q\left(  2e_{i}\right)
}.\label{TJMaxx}%
\end{equation}
Thus, recalling condition (\ref{eq:cdt-q-discrete}), and chosing $t=C\log
^{2}\beta/\beta^{2}$, we obtain
\begin{align}
\frac{1}{t}\mathbf{E}(\log Z_{t}) &  \geq\frac{\log t}{t}-2d+\frac{2\beta
}{3\sqrt{\pi t}}c_{Q}\nonumber\\
&  \geq\frac{\beta^{2}}{\log\beta}\left(  -\frac{2}{C}+\frac{2c_{Q}}%
{3\sqrt{c\pi}}\right)  +\frac{\beta^{2}}{C\log^{2}\beta}\left(  \log
C+2\log\log\beta\right)  -2d.\label{eq:low-bnd-zt-discrete-interm}%
\end{align}
The proof is completed by choosing $C$ such that $-\frac{2}{C}+\frac{2c_{Q}%
}{3\sqrt{C\pi}}>0$, i.e. $C>\frac{9\pi}{Q(0)-Q(2)}$, and $\beta$ large enough
so that the second and third terms in (\ref{eq:low-bnd-zt-discrete-interm})
contribute nonnegatively.
\end{proof}

\subsection{Upper bound for the random walk polymer}

The upper bound result in Theorem \ref{thm:bnd-free-nrj-walk} can be
summarized in the following proposition.

\begin{proposition}
Under the assumption that $Q\left(  0\right)  <\infty$, there exists a
constant $\beta_{0}^{\prime}>0$, which depends on $Q$ and on $d$, and a
constant $C>0$, which depend only on $Q$, such that if $\beta>\beta
_{0}^{\prime}$ then almost surely
\[
\lim_{t\rightarrow\infty}\frac{1}{t}\log Z_{t}\leq Cd^{3}\frac{\beta^{2}}%
{\log\beta}.
\]

\end{proposition}

\begin{proof}
Define $\mathcal{S}(t,n)$, $k_{j}$, $t_{i}^{j}$, $\tilde{t}_{l}$, $N_{t}$,
$\mathcal{P}_{n}$ and $\tilde{x}_{l}$ like in Step 3 of the proof of
Proposition \ref{prop:upper-bnd-holder}. Then if we fix $N_{t}=m$, we can
define
\[
X\left(  m,\tilde{t},\tilde{x}\right)  :=\sum_{i=0}^{m}\left\{  W\left(
\tilde{t}_{i+1},\tilde{x}_{i}\right)  -W\left(  \tilde{t}_{i},\tilde{x}%
_{i}\right)  \right\}  .
\]
Let $\alpha$ be a fixed positive number which will be chosen later. Let
$I_{\alpha}=\cup_{m\leq\alpha t}J_{m}$, where $J_{m}:=\left\{  m\right\}
\times\mathcal{S}_{m,t}\times\mathcal{P}_{m}$, and set also $Y_{\alpha}%
=\sup_{I_{\alpha}}X.$ As in the Brownian case, we can bound $\mathbf{E}\left[
\log Z_{t}\right]  $ above as follows:
\begin{equation}
\mathbf{E}\left[  \log Z_{t}\right]  \leq\mathbf{E}\left[  \log\left(
A+B\right)  \right]  \leq\mathbf{E}\left[  \log_{+}A\right]  +\mathbf{E}%
\left[  \log_{+}B\right]  +\log2,\label{Uub}%
\end{equation}
where $\log_{+}A=\left(  \log A\right)  _{+}=\max(\log A,0)$ and
\begin{align}
A &  :=P_{b}\left[  N_{t}\leq\alpha t\right]  \exp\left(  \beta Y_{\alpha
}\right)  \label{eq:def-A-B-discrete}\\
B &  :=\sum_{m>\alpha t}P_{b}\left[  N_{t}=m\right]  \,E_{b}\left[
\exp\left(  \beta X\left(  m,\tilde{t},\tilde{x}\right)  \right)
\Big|\Big.N_{t}=m\right]  .\nonumber
\end{align}

\noindent\textit{Step 1: The term A.} As in the continuous case, we have that
\begin{equation}
\mathbf{E}\left[  \sup_{T_{n\alpha}}X(k,\tilde{t},\tilde{x})\right]  \leq
Ktd\sqrt{n\alpha},\label{stimasup2}%
\end{equation}
where $K$ depends only on $Q$. So, bounding $P_{b}\left[  N_{t}\leq\alpha
t\right]  $ by $1$, we have
\begin{equation}
\mathbf{E}\left[  \log_{+}A\right]  \leq\beta\mathbf{E}\left[  Y_{\alpha
}\right]  \leq\beta Kdt\sqrt{\alpha}.\label{A2}%
\end{equation}

\noindent\textit{Step 2: The term B.} The term $B$ defined in
(\ref{eq:def-A-B-discrete}) can be bounded as follows:
\begin{align*}
&  \mathbf{E}\left[  \log B_{+}\right] \\
&  =\mathbf{E}\left[  \log_{+}\sum_{m>\alpha t}P_{b}\left[  N_{t}=m\right]
\sum_{\tilde{x}\in\mathcal{P}_{m}}\frac{1}{(2d)^{m}}\int_{\mathcal{S}_{m,t}%
}\exp\left(  \beta X\left(  m,\tilde{t},\tilde{x}\right)  \right)  d\tilde
{t}\right] \\
&  =\mathbf{E}\left[  \log_{+}\sum_{n\geq1}\sum_{m\in\lbrack\alpha
nt,\alpha(n+1)t]}P_{b}\left[  N_{t}=m\right]  \sum_{\tilde{x}\in
\mathcal{P}_{m}}\frac{1}{(2d)^{m}}\int_{\mathcal{S}_{m,t}}\exp\left(  \beta
X\left(  m,\tilde{t},\tilde{x}\right)  \right)  d\tilde{t}\right] \\
&  \leq\mathbf{E}\left[  \log_{+}\sum_{n\geq1}P_{b}\left[  N_{t}>\alpha
nt\right]  \exp\left(  \beta Y_{(n+1)\alpha}\right)  \right]  .
\end{align*}
So, using the fact that for $t>1$, the power $t^{-1}$ of a sum is less than
the sum of the terms raised to the power $t^{-1}$, followed by Jensen's
inequality, we have, similarly to what we did in the proof of Proposition
\ref{prop:upper-bnd-holder},
\[
\frac{1}{t}\mathbf{E}\left[  \log_{+}B\right]  \le\log\left(  1+\sum_{n\geq
1}\left(  P_{b}\left[  N_{t}>\alpha nt\right]  \right)  ^{t^{-1}}%
\mathbf{E}\left[  \exp\left(  \frac{\beta Y_{(n+1)\alpha}}{t}\right)  \right]
\right)  .
\]

Using once again Gaussian supremum analysis results (see \cite{A} or \cite{VV}), for any
$\alpha,x>0$,
\begin{align*}
\mathbf{E}\left[  \exp\left(  xY_{\alpha}\right)  \right]   &  \leq\exp\left(
x\mathbf{E}\left[  Y_{\alpha}\right]  \right)  \exp\left(  x^{2}K_{u}%
\max_{(m,\tilde{t},\tilde{x})\in I_{\alpha}}\mbox{\bf Var}\left[  X\left(
m,\tilde{t},\tilde{x}\right)  \right]  \right)  \\
&  \leq\exp\left(  xdKt\sqrt{\alpha}\right)  \exp\left(  x^{2}K_{u}%
tQ(0)\right),
\end{align*}
where $K_u$ designate a universal constant, and
where we used (\ref{stimasup2}) and the trivial fact $\mathbf{E}[X\left(
m,\tilde{t},\tilde{x}\right)  ^{2}]=Q\left(  0\right)  t$. Hence%
\[
\mathbf{E}\left[  \exp\left(  \frac{\beta Y_{(n+1)\alpha}}{t}\right)  \right]
\leq\exp\left(  \beta dK\sqrt{\alpha(n+1)}+\frac{\beta^{2}K_{u}Q(0)}%
{t}\right)  .
\]
If we choose $t$ such that $t>(2\beta K_{u}Q(0))/(dK\alpha^{1/2})$, the
estimate on $B$ becomes
\begin{equation}
\frac{1}{t}\mathbf{E}\left[  \log_{+}B\right]  \leq\log\left\{  1+\sum
_{n\geq1}\left(  P_{b}\left[  N_{t}>\alpha nt\right]  \right)  ^{t^{-1}}%
\exp\left(  \beta dK_{u}\sqrt{\alpha}\left(  \sqrt{n+1}+\frac{1}{2}\right)
\right)  \right\}  .\label{logB+}%
\end{equation}

\noindent\textit{Step 3: The tail of }$N_{t}$\textit{.} Using the presumably
well-known tail estimate $P_{b}\left[  N_{t}>\alpha t\right]  $ $\leq
\exp\left(  -\alpha t\log\left(  \frac{\alpha}{2d}\right)  -t(\alpha
-2d)\right)  $, valid for all $\alpha\geq1$ (see e.g. \cite[pages
16-19]{Massart}), if we set $\alpha^{\prime}=\alpha/2d$ and we assume
$\alpha^{\prime}\geq\exp\left(  1-1/2d\right)  $ we have
\begin{equation}
P_{b}\left[  N_{t}>\alpha t\right]  \leq\exp\left(  -t\alpha^{\prime}%
\log\alpha^{\prime}\right)  .\label{tailestim}%
\end{equation}

\noindent\textit{Step 4: Grouping our estimates and choosing }$\alpha
$\textit{.} From (\ref{logB+}) and (\ref{tailestim}) we have
\[
\frac{1}{t}\mathbf{E}\left[  \log_{+}B\right]  \leq\log\left\{  1+\sum
_{n\geq1}\exp\left(  -\alpha^{\prime}n\log\alpha^{\prime}n+d\beta
K_{u}\sqrt{\alpha}\left(  \sqrt{n+1}+\frac{1}{2}\right)  \right)  \right\}
.
\]
To exploit the negativity of the exponential term, we simply require%
\begin{equation}
\alpha^{\prime}\log\alpha^{\prime}=4dK_{u}\beta\sqrt{\alpha}%
.\label{alphadef}%
\end{equation}
Indeed, since $n\geq1$, we then have that the term inside the exponential is
\begin{align*}
&  \alpha^{\prime}n\log\alpha^{\prime}n-\beta dK_{u}\sqrt{\alpha}\left(
\sqrt{n+1}+\frac{1}{2}\right)  =\alpha^{\prime}n\log\alpha^{\prime}n-\frac
{1}{4}\alpha^{\prime}\log\alpha^{\prime}\left(  \sqrt{n+1}+\frac{1}{2}\right)
\\
&  \geq\alpha^{\prime}n\log\alpha^{\prime}-\frac{1}{4}\alpha^{\prime}%
\log\alpha^{\prime}\left(  \sqrt{n+1}+\frac{1}{2}\right)  \\
&=\frac{1}{2}%
\alpha^{\prime}n\log\alpha^{\prime}n+\left(  \frac{1}{2}n-\frac{1}{4}\left(
\sqrt{n+1}+\frac{1}{2}\right)  \right)  \alpha^{\prime}\log\alpha^{\prime}\\
&  \geq\frac{1}{2}\alpha^{\prime}n\log\alpha^{\prime},
\end{align*}
which implies
\[
\frac{1}{t}\mathbf{E}\left[  \log_{+}B\right]  \leq\log\left\{  1+\sum
_{n\geq1}\exp\left(  -\frac{1}{2}\alpha^{\prime}n\log\alpha^{\prime}\right)
\right\}  =\log\left\{  1+\frac{1}{\exp\left(  \frac{1}{2}\alpha^{\prime}%
\log\alpha^{\prime}\right)  -1}\right\}  :=c_{d}.
\]
The restriction $\alpha^{\prime}\geq\exp(1-1/2d)$ implies that $c_{d}$ is a
constant that depends on the dimension $d$ only. Combining this with
(\ref{Uub}) and (\ref{A2}), we get
\begin{equation}
\frac{1}{t}\mathbf{E}\left[  \log Z_{t}\right]  \leq\frac{\log2}{t}%
+c_{d}+dK_{u}\beta\sqrt{\alpha}.\label{sayepresk}%
\end{equation}

\noindent\textit{Step 5: Conclusion.} It is easy enough to see that, with
\begin{equation}
x:=\left(  4d\sqrt{2d}\beta K_{u}\right)  ^{2},\label{xdef}%
\end{equation}
the equation (\ref{alphadef}), which is $\alpha^{\prime}=x/\log^{2}%
\alpha^{\prime}$, has a unique solution $\alpha^{\prime}$ when $x$ exceeds
$e$, and that $\alpha^{\prime}$ also exceeds $e$ in that case: indeed
$\alpha^{\prime}=e$ when $x=e$ and $d\alpha^{\prime}/dx=\left(  \log^{2}%
\alpha^{\prime}+2\log\alpha^{\prime}\right)  ^{-1}>0$ for all $\alpha^{\prime
}\geq e$. Therefore, since $\log^{2}\alpha^{\prime}>1$, we can write
$\alpha^{\prime}\leq x$, and thus we also have:
\begin{equation}
\alpha^{\prime}=\frac{x}{\log^{2}\alpha^{\prime}}\geq\frac{x}{\log^{2}%
x}.\label{alphalb}%
\end{equation}
This lower bound on $\alpha^{\prime}$ implies the following upper bound on
$\alpha^{\prime}$:%
\begin{equation}
\alpha^{\prime}=\frac{x}{\log^{2}\alpha^{\prime}}\leq\frac{x}{\log^{2}\left(
x/\log^{2}x\right)  }=\frac{x}{\left(  \log x-2\log(\log x)\right)  ^{2}%
}.\label{eq:exp-alpha-prime}%
\end{equation}
Since there exists $x_{0}$ such that, for any $x>x_{0}$, we have
\begin{equation}
\log x>4\log(\log x),\label{eq:bnd2-beta-discrete}%
\end{equation}
and we can recast expression (\ref{eq:exp-alpha-prime}) into:
\[
\alpha^{\prime}\leq\frac{4x}{\log^{2}x}=\frac{(4d\sqrt{2d}\beta K_{u})^{2}%
}{(\log\beta+\log4d\sqrt{2d}K_{u})^{2}}\leq(4d\sqrt{2d}K_{u})^{2}%
\frac{\beta^{2}}{\log^{2}\beta},
\]
from which we obtain $\alpha\leq(8d^{2}K_{u})^{2}\beta^{2}\log^{-2}\beta$ (recall
$\alpha=:2d\alpha^{\prime}$). Thus for $t$ large enough:
\[
\frac{\mbox{{\bf E}}[\log(Z_{t})]}{t}\leq\frac{\log2}{t}+c_{d}+8d^{3}%
K_{u}^{2}\frac{\beta^{2}}{\log\beta}.
\]
Taking limits as $t$ tends to $\infty$ and choosing $\beta$ so that
\begin{equation}
\beta^{2}>\frac{c_{d}\log\beta}{8d^{3}K_{u}^{2}}%
,\label{eq:bnd3-beta-discrete}%
\end{equation}
the theorem is proved with $C=16K_{u}^{2}$.\vspace{0.3cm}

Finally, we show that the theorem holds for $\beta$ large enough. Analyzing
the conditions we used above, we only have to take $\beta\geq\beta_{0}%
^{\prime}:=\max(K_{d},\beta^{\ast},\beta_{\ast}),$where $K_{d},\beta^{\ast
},\beta_{\ast}$ are now specified. This is due to the fact that we assumed
$\alpha^{\prime}\geq\exp(1-1/2d)$ and this implies, via (\ref{alphalb}), that
$x\geq2d\exp(1-1/2d)\left(  \log2d+1-1/\left(  2d\right)  \right)  $, and
therefore, from (\ref{xdef}), we have to take $\beta\geq K_{d}$, where $K_{d}$
is a constant that depends only on the dimension $d$. In addition, according
to (\ref{eq:bnd2-beta-discrete}) and (\ref{eq:bnd3-beta-discrete}),
$\beta_{\ast}$ and $\beta^{\ast}$ are the solutions to the following
equations:
\[
\beta^{2}=\frac{c_{d}\log\beta}{8d^{3}K_{u}^{2}}\quad\mbox{ and }\quad
\log(4d\sqrt{2d}K_{u}\beta)=4\log(\log(4d\sqrt{2d}K_{u}\beta)).
\]

\end{proof}

\end{document}